\documentclass[a4paper,11pt,twoside]{article}

\usepackage{xypic}
\usepackage{latexsym}
\usepackage{amsmath}
\usepackage{amstext}
\usepackage{amsfonts}
\usepackage{amssymb}
\usepackage{amsbsy}
\usepackage{enumerate}
\usepackage{theorem}

\xyoption{all}
\CompileMatrices

\newtheorem{theorem}{Theorem}

\newtheorem{lemma}{Lemma}[section]
\newtheorem{lemma-def}[lemma]{Lemma and Definition}
\newtheorem{proposition}[lemma]{Proposition}
\newtheorem{corollary}[lemma]{Corollary}

\theorembodyfont{\rmfamily}
\newtheorem{remark}[lemma]{Remark}
\newtheorem{definition}[lemma]{Definition}

\newcommand{\id}{\text{\rm id}}

\newcommand{\Aut}{\text{\rm Aut}}

\newcommand{\Iso}{\text{\rm Iso}}

\newcommand{\Hom}{\text{\rm Hom}}

\newcommand{\Ann}{\text{\rm Ann}}

\newcommand{\Spec}{\text{\rm Spec}}
\newcommand{\Graph}{\text{\rm Graph}}

\renewcommand{\Im}{\text{\rm Im}}

\newcommand{\qed} {\hbox{} \nolinebreak \hfill $\;\Box$}

\begin{document}

\title{Families of elliptic curves with genus 2 covers of degree 2}

\author{Claus Diem}
\date{November 8, 2005}

\maketitle

\pagestyle{myheadings}
\markboth{\sc Diem}{\sc Elliptic curves with genus 2 covers of degree 2}

\begin{abstract}
We study genus 2 covers of relative elliptic curves over an arbitrary base in which 2 is invertible. Particular emphasis lies on the case that the covering degree is 2. We show that the data in the "basic construction" of genus 2 covers of relative elliptic curves determine the cover in a unique way (up to isomorphism). 

A classical theorem says that a genus 2 cover of an elliptic curve of degree 2 over a field of characteristic $\neq 2$ is birational to a product of two elliptic curves over the projective line. We formulate and prove a generalization of this theorem for the relative situation.

We also prove a Torelli theorem for genus 2 curves over an arbitrary base.
\end{abstract}

\vspace{1 ex}
\begin{center}\begin{minipage}{110mm}\footnotesize{\bf Key words:} Elliptic curves, covers of curves, families of curves, curves of genus 2, curves with split Jacobian.
\end{minipage}
\end{center}

\begin{center}\begin{minipage}{110mm}\footnotesize{\bf MSC2000:} 14H45, 14H10, 14H30.

\end{minipage}
\end{center}
\vspace{1 ex}

\section*{Introduction}
The purpose of this article is to study covers $f : C \longrightarrow E$ where $C/S$ is a (relative, smooth, proper) genus 2 curve, $E/S$ is a (relative) elliptic curve and the base $S$ is a locally noetherian scheme over $\mathbb{Z}[1/2]$. Particular emphasis lies on the case that the covering degree $N$ is 2.

If one studies genus 2 covers of (relative) elliptic curves, it is convenient to restrict ones attention to so-called \emph{minimal} covers. These are covers $C \longrightarrow E$ which do not factor over a non-trivial isogeny $\tilde{E} \longrightarrow E$. If now $f : C \longrightarrow E$ is a minimal cover and $x \in E(S)$, then $T_x \circ f$ is also one. This ambiguity motivates the notion of a \emph{normalized} cover introduced in \cite{KaHur}: By definition, such a cover is minimal and satisfies a certain condition concerning the direct image of the Weierstra{\ss} divisor of $C$ on $E$ (for precise definition see below). Now for every minimal cover $f : C \longrightarrow E$ there is exactly one $x \in E(S)$ such that $T_x \circ f : C \longrightarrow E$ is normalized.

To every minimal cover $f : C \longrightarrow E$ one can associate in a canonical way an elliptic curve $E'_f/S$ and an isomorphism of $S$-group schemes $\psi_f : E[N] \tilde{\longrightarrow} E'_f[N]$ which is anti-isometric with respect to the Weil pairing; see \cite{KaHur}. It is shown in \cite{KaHur} that for fixed $S$, $E/S$ and $N \geq 3$, the assignment $f \mapsto (E_f,\psi_f)$ induces a monomorphism from the set of isomorphism classes of normalized genus 2 covers of degree $N$ of $E/S$ to the set of isomorphism classes of tuples $(E', \psi)$ of elliptic curves $E/S$ with an anti-isometric isomorphism $\psi : E[N] \tilde{\longrightarrow} E'[N]$. Explicit conditions are given when a tuple $(E',\psi)$ corresponds to a normalized genus 2 cover $C \longrightarrow E$ of degree $N$ over $S$ -- this is called ``basic construction'' in \cite{KaHur}.

In this work, we show that the above assignment is in fact a monomorphism for all $N \geq 2$. Our starting point is a Torelli theorem (Theorem \ref{genus-2-torelli}) for relative genus 2 curves which follows rather easily from the detailed appendix of \cite{KaHur}. With the help of this theorem, we prove a Torelli theorem for normalized genus 2 covers of (relative) elliptic curves; see Proposition \ref{genus-2-cover-Torelli}. This result implies immediately that the ``Torelli map'' of \cite{KaHur} is a monomorphism for arbitrary $N \geq 2$. In \cite{KaHur}, the corresponding statement is only proved for $N\geq 3$ and the proof is more involved; cf. \cite[Proposition 5.12]{KaHur}. The injectivity of the above assignment then follows with other results of \cite{KaHur}.

For $N=2$ (and fixed $S$ and $E/S$), tuples $(E',\psi)$ as well as normalized covers $C \longrightarrow E$ have a non-trivial automorphism of order 2. This leads to a certain ``non-rigidity'' in the ``basic construction'': Any two covers corresponding to the same tuple $(E',\psi)$ are isomorphic, but the isomorphism is not unique. We propose a ``symmetric basic construction'' which leads to a more rigid statement (and is more explicit than the ``basic construction'').

We then fully concentrate on the case that $N=2$. We show in particular that for every normalized cover $f: C \longrightarrow E$ of degree 2, one has a canonical commutative diagram
\[ \xymatrix{
& \ar_f[dl] C \ar[dr] & \\
E \ar[dr] & & \ar[dl] E'_f \; , \\
& {\mathbf{P}} &}
\]
where $\mathbf{P} := E/\langle [-1] \rangle$ is a $\mathbb{P}^1$-bundle over $S$ and all morphisms are covers of degree 2 such that the induced morphism $C \longrightarrow E \times_{\mathbf{P}} E'_f$ induces birational morphisms on the fibers over $S$; see Theorem \ref{theorem-P-diagram} in Section \ref{section:genus2} and Corollary \ref{corollary-P-diagram}. This generalizes a classical result on genus 2 curves with elliptic differentials of degree 2 over a field of characteristic $\neq 2$ which follows immediately from Kummer theory applied to the extension $\kappa(C)/\kappa(E/\langle [-1] \rangle )$.

Finally, we discuss a reinterpretation of this result and show that it is closely related to a general statement on $\mathbb{P}^1$-bundles which we prove in an appendix.\\

The study of genus 2 curves with split Jacobian has a long history which arguably started with the task of reducing hyperelliptic integrals of genus 2 of the first kind to sums of elliptic integrals. Here a substitution of variables gives rise to a genus 2 cover of an elliptic curve. The study for degree 2 dates back to Legendre who gave the first examples and Jacobi. More information on this classical material can be found in \cite{Kra}, pp.477-482.

It is now also classical that to every minimal cover $f : C \longrightarrow E$ one can in a canonical way associate a ``complementary'' minimal cover $C \longrightarrow E'_f$ of the same degree (unique up to translation on $E$); see e.g.\ \cite{Kuhn}. The idea to describe genus 2 covers of a fixed elliptic curve $E$ (over a field) by giving the complementary elliptic curve $E'_f$ and a suitable anti-isometric isomorphism $E[N] \tilde{\longrightarrow} E'_f[N]$, where $N$ is the covering degree, is due to G.\ Frey and E.\ Kani; see \cite{FKBasicC} and also \cite{KaNum}. The basic results for genus 2 covers of \emph{relative} elliptic curves were obtained by E.\ Kani in \cite{KaHur}.\\

An application of some results presented in this article can be found in \cite{DF}. In this work, examples of relative, non-isotrivial genus 2 curves $C/S$ which possess an infinite tower of non-trivial \'etale covers $\cdots \longrightarrow C_i \longrightarrow \ldots C_0 = C$ such that for all $i$, $C_i \longrightarrow C$ is Galois and $C_i/S$ is also a curve (in particular has \emph{geometrically connected} fibers) are given. The genus 2 curves in question are covers of elliptic curves with covering degree 2, the base schemes are affine curves over finite fields of odd \nolinebreak characteristic.

\subsection*{Terminology and notation}
This work is closely related to \cite{KaHur}. With the exception of the following assumption, the following three definitions and Definition \ref{def:symmetric-pair}, all definitions and notations follow this work. We thus advise the reader to have \cite{KaHur} at hand when he goes through the details of this article. Note that although the primary emphasis of \cite{KaHur} lies on genus 2 covers of elliptic curves $E_S$, where $E/K$ is an elliptic curve over a field $K$ of characteristic $\neq 2$ and $S$ is a $K$-scheme, as stated in various places of \cite{KaHur}, the results of \cite{KaHur} hold for genus 2 covers of elliptic curves over arbitrary locally noetherian schemes over $\mathbb{Z}[1/2]$.

If not stated otherwise, all schemes we consider are assumed to be locally noetherian.

If $g \in \mathbb{N}_0$, then a \emph{(relative) curve} of genus $g$ over $S$ is a smooth, proper morphism $C \longrightarrow S$ whose fibers are geometrically connected curves of genus $g$. (We thus do not assume that the genus is $\geq 1$ or that for $g=1$ $C/S$ has a section.)

If $C/S$ is a curve and $N \in \mathbb{N}$, $g \in \mathbb{N}_0$, then a genus $g$ \emph{cover} of degree $N$ of $C$ is an $S$-morphism $f : C' \longrightarrow C$, where $C'/S$ is a genus $g$ curve, which induces morphisms of the same degree $N$ on the fibers over $S$. (Note that $f$ is automatically finite, flat and surjective; cf.\ \cite[Section 7, 7)]{KaHur}.)

If $C/S$ and $C'/S$ are two curves of genus $\geq 2$, we denote the scheme of $S$-isomorphisms from $C$ to $C'$ by $\mathbf{Iso}_S(C,C')$; cf.\ \cite{DM}.

Following \cite{LK}, a curve $C/S$ is called \emph{hyperelliptic} if it has a (by Lemma \ref{red-inj} necessarily unique) automorphism $\sigma_{C/S}$ which induces hyperelliptic involutions on the geometric fibers. For equivalent definitions of $\sigma_{C/S}$, see \cite[Theorem 5.5]{LK}.

We have used the following definition in the introduction; cf.\ \cite{KaHur}:

Let $S$ be a scheme over $\mathbb{Z}[1/2]$, let $C/S$ be a genus 2 curve and let $E/S$ be an elliptic curve. Then a cover $f : C \longrightarrow E$ is \emph{minimal} if it does not factor over a non-trivial isogeny $\tilde{E} \longrightarrow E$, and it is \emph{normalized} if it is minimal and we have the equality of relative effective Cartier divisors
\[ f_*(W_{C/S}) = 3 \epsilon [0_{E/S}] + (2- \epsilon) E[2]^{\#}\; ,\]
where $W_{C/S}$ is the Weierstra{\ss} divisor of $C/S$, $E[2]^{\#} := E[2] - [0_{E/S}]$ and $\epsilon = 0$ if $\deg(f)$ is even and $\epsilon = 1$ if $\deg(f)$ is odd.\footnote{There are misprints in the definitions in \cite[Section 2]{KaHur} and \cite[Section 3]{KaHur}.} Note that a normalized cover satisfies
\begin{equation}
\label{pseudo-normalized}
f \circ \sigma_{C/S} = [-1] \circ f \; ;
\end{equation}
cf.\ \cite[Theorem 3.2 (c)]{KaHur}.

We frequently use the following notation:

If $f: T \longrightarrow S$ is a morphism of schemes and $\varphi: X \longrightarrow Y$ is a morphism of $S$-schemes, we denote the morphism induced by base change via $f$ by $f^* \varphi : f^* X \longrightarrow f^* Y$ or just $\varphi_{T} : X_{T} \longrightarrow Y_{T}$.

We use two different symbols to denote isomorphisms: If we just want to state that two objects $X,Y$ in some category are isomorphic, we write $X \approx Y$. If $X$ and $Y$ are isomorphic with respect to a canonical isomorphism or with respect to a fixed isomorphism which is obvious from the context, we write $X \simeq Y$.\\

\paragraph{Acknowledgments.}
The author would like to thank G.\ Frey, E.\ Kani and E.\ Viehweg for various discussions related to this work.

\section{A Torelli theorem for relative genus 2 curves}

The purpose of this section is to prove the following theorem.

\begin{theorem}
\label{genus-2-torelli}
Let $S$ be a scheme, let $C/S$ and $C'/S$ be two genus 2 curves. Then the map $\Iso_S(C, C') \longrightarrow \Iso_S((J_{C}, \lambda_{C}), (J_{C'},\lambda_{C'})), \; \varphi \mapsto \varphi_*$
is an isomorphism.
\end{theorem}
Here, by $\lambda_{C}$ we denote the canonical polarization of the Jacobian $J_C$ of a genus 2 curve $C/S$ and for an isomorphism $\varphi : C \longrightarrow C'$ of two genus 2 curves over $S$, we define $\varphi_* := \lambda_{C'}^{-1} \circ (\varphi^*\hat{)} \circ \lambda_{C} = (\varphi^*)^{-1}$.

This Torelli theorem for (relative) genus 2 curves is well known in the case that $S$ is the spectrum of an (algebraically closed) field; cf.\ e.g.\ \cite[Theorem 12.1]{Mi-JV} where it is stated with a slightly different formulation for arbitrary hyperelliptic curves over algebraically closed fields.

Theorem \ref{genus-2-torelli} follows from Lemmata \ref{torelli-injective} and \ref{torelli-surjective} which are proved below.

Let $S$ be a scheme, and let $C/S$ and $C'/S$ be curves.

We will frequently use the fact that the formation of the Jacobian commutes with arbitrary base-change: Let $f: T \longrightarrow S$ be a morphism of schemes. Then we have canonical isomorphisms $(J_{C_T},\lambda_{C_T}) \simeq ((J_C)_T,(\lambda_{C})_T)$, $(J_{C'_T},\lambda_{C'_T}) \simeq ((J_{C'})_T,(\lambda_{C'})_T)$. Moreover, under the obvious identification, we have
\begin{equation}
\label{pull-back-equality}
(\varphi_*)_T = (\varphi_T)_* : J_{C_T} \longrightarrow J_{C'_T} \; \; \text{i.e.} \; \; f^*(\varphi_*) = (f^* \varphi)_* \; .
\end{equation}

\begin{lemma}
\label{red-inj}
Let $S$ be a connected scheme, let $s \in S$. Then the restriction map $\Iso_S(C,C') \longrightarrow \Iso_{\kappa(s)}(C_s,C'_s)$ is injective.
\end{lemma}
\emph{Proof.}
The $S$-isomorphisms between $C$ and $C'$ correspond to sections of the $S$-scheme $\textbf{Iso}_S(C,C')$. As this scheme is unramified over $S$ (see \cite[Theorem 1.11]{DM}), the result follows with \cite[Expos\'e I, Corollaire 5.3.]{SGA}.
\qed

\begin{lemma}
\label{torelli-injective}
Let $S$ be a connected scheme, let $s \in S$. Then the map $\Iso_S(C, C') \longrightarrow \Iso_{\kappa(s)}((J_{C_s}, \lambda_{C_s}), (J_{C'_s},\lambda_{C'_s})), \; \varphi \mapsto (\varphi_s)_* = (\varphi_*)_s$ is injective.
\end{lemma}
\emph{Proof.}
This follows from the previous lemma and the classical Torelli Theorem (see \cite[Theorem 12.1]{Mi-JV}).\qed

\begin{lemma}
\label{descent}
Let $S' \longrightarrow S$ be faithfully flat and quasi compact. Let $\varphi' : C_{S'} \longrightarrow C'_{S'}$ be an $S'$-isomorphism, and let $\alpha : J_C \longrightarrow J_{C'}$ be a homomorphism with $\alpha_{S'} = \varphi'_*$. Then there exists an $S$-isomorphism $\varphi : C \longrightarrow C'$ with $\varphi_{S'} = \varphi'$ and $\alpha = \varphi_*$.
\end{lemma}
\emph{Proof.}
Let $S'' := S' \times_S S'$, let $p_1, p_2 : S'' \longrightarrow S'$ be the two projections. We want to show that $p_1^*\varphi' = p_2^*\varphi'$. Then the statement follows by faithfully flat descent; see \cite[Section 6.1., Theorem 6]{BLR}.

By assumption we have $p_1^*(\varphi'_*) = p_2^*(\varphi'_*)$. Together with (\ref{pull-back-equality}) this implies that $(p_1^* \, \varphi')_* = (p_2^* \, \varphi')_*$. Now the equality $p_1^*\varphi' = p_2^*\varphi'$ follows with the previous lemma.
\qed

The following lemma is a special case of \cite[Proposition 6.1]{Mu-GIT}, the ``Rigidity Lemma''.

\begin{lemma}
\label{rigidity-lemma-special-case}
Let $S$ be a connected scheme, let $s \in S$. Let $A/S, A'/S$ be two abelian schemes. Then the map
$\Hom_S(A,A') \longrightarrow \Hom_{\kappa(s)}(A_s,A'_s)$ is injective.
\end{lemma}

\begin{lemma}
Let $C/S$ and $C'/S$ be genus 2 curves, and assume that both curves have a section. Then the map $\Iso_S(C, C') \longrightarrow \Iso_S((J_{C}, \lambda_{C}), (J_{C'},\lambda_{C'})), $\linebreak$\varphi \mapsto \varphi_*$ is surjective.
\end{lemma}
\emph{Proof.}
Let $a : S \longrightarrow C$ be a section. Let $j_a : C \longrightarrow J_C$ be the immersion associated to $a$; cf. \cite[Section 7, 6)]{KaHur}. Analogously, let $a' : S \longrightarrow C'$ be a section, and let $j_{a'} : C' \longrightarrow J_{C'}$ be the associated immersion. Now $j_a(C)$ is a Cartier divisor on $J_C$ which defines the principal polarization $\lambda_{C}$. (Indeed, for all $s \in S$, we have $\lambda_{C_s} = \lambda_{\mathcal{O}(j_a(C)_s)} : J_{C_s} \longrightarrow J_{C'_s}$. The equality $\lambda_{C} = \lambda_{\mathcal{O}(j_a(C))}$ follows with Lemma \ref{rigidity-lemma-special-case}.) Analogously, $j_{a'}(C')$ is an a Cartier divisor on $J_{C'}$ which defines the principal polarization $\lambda_{C'}$.

Let $\alpha : J_C \longrightarrow J_{C'}$ be an isomorphism which preserves the principal polarizations, i.e.\ which satisfies $\hat{\alpha} \circ \lambda_{C'} \circ \alpha = \lambda_{C}$.

Then $\lambda_{C}$ is given by the divisor $\alpha^{-1}(j_{a'}(C'))$. It follows from \cite[Lemma 7.1]{KaHur} that $\alpha^{-1}(j_{a'}(C')) = T_x^{-1}(j_a(C))$ for some $x \in J_C(S)$. This can be rewritten as $(\alpha^{-1} \circ j_{a'})(C') = (T_{-x} \circ j_a)(C)$. Note here that $\alpha^{-1} \circ j_{a'} : C' \longrightarrow J_C$ and $T_{-x} \circ j_a : C \longrightarrow J_C$ are closed immersions, and we have an equality of the associated closed subschemes of $J_{C'}$. This means that there exists an isomorphism of schemes $\varphi : C \longrightarrow C'$ such that $\alpha^{-1} \circ j_{a'} \circ \varphi = T_{-x} \circ j_a$, i.e.\ $j_{a'} \circ \varphi = \alpha \circ T_{-x} \circ j_a$. A short calculation shows that $\varphi$ is in fact an $S$-isomorphism.

The equality $j_{a'} \circ \varphi = \alpha \circ T_{-x} \circ j_a$ immediately implies that $\varphi_* = \alpha$.
\qed

\begin{lemma}
\label{torelli-surjective}
Let $C/S$, $C'/S$ be two genus 2 curves. Then the map \linebreak $\Iso_S(C, C') \longrightarrow \Iso_S((J_{C}, \lambda_{C}), (J_{C'},\lambda_{C'})), \; \varphi \mapsto \varphi_*$ is surjective.
\end{lemma}
\emph{Proof.}
Let $W_{C/S}$, $W_{C'/S}$ be the Weierstra{\ss} divisors of $C/S$ and $C'/S$ respectively and let $W:= W_{C/S} \times_S W_{C'/S}$. Now the canonical map $W \longrightarrow S$ is faithfully flat and quasi compact (in fact it is finite flat of degree 36), and $C_W/W$ as well as $C'_W/W$ have sections (namely the sections induced by $W_{C/S} \hookrightarrow C$, $W_{C'/S} \hookrightarrow C'$). It follows by the above lemma that $\Iso_{W}(C_{W}, C'_{W}) \longrightarrow \Iso_{W}((J_{C_{W}}, \lambda_{C_{W}}), (J_{C_W},\lambda_{C_W})), \, \varphi \mapsto \varphi_*$
is surjective. The claim now follows with Lemma \ref{descent}.
\qed\\

The above considerations easily imply:

\begin{corollary}
\label{hyerelliptic-involution-commutativity}
Let $C/S, C'/S$ be hyperelliptic curves, let $\varphi : C \longrightarrow C'$ be an $S$-isomorphism. Then
\[ \sigma_{C'/S} \circ \varphi = \varphi \circ \sigma_{C/S}\; . \]
\end{corollary}
\emph{Proof.}
We can assume that $S$ is connected. Let $s \in S$. It is well known that $(\sigma_{C_s})_* = [-1], (\sigma_{C'_s})_* = [-1]$. This implies $(\sigma_{C'_s})_* \circ (\varphi_s)_* = - (\varphi_s)_* = (\varphi_s)_* \circ (\sigma_{C_s})_*$. The result now follows with Lemma \ref{torelli-injective}.
\qed

\smallskip

We also have:

\begin{lemma}
\label{[-1]}
Let $C/S$ be a hyperelliptic curve. Then $(\sigma_{C/S})_* = [-1]$.
\end{lemma}
\emph{Proof.}
This follows from the well known result over the spectrum of a field by Lemma \ref{rigidity-lemma-special-case}.
\qed

\section{Review of the ``basic construction''}

Theorem \ref{genus-2-torelli} can be used to prove a Torelli theorem for normalized genus 2 covers of elliptic curves which in turn can be used to simplify some proofs in \cite{KaHur} as well as to strengthen the results for the case that the covering degree $N$ is $2$. This is done in the first half of this section. Throughout the section, we freely use results from \cite{KaHur}.

Let $S$ be a scheme over $\mathbb{Z}[1/2]$. The following definition is analogous to the ``notation'' in Section 3 of \cite{KaHur}.

\begin{definition}
Let $E/S$ be an elliptic curve. Let $f_1 : C_1 \longrightarrow E, f_2 : C_2 \longrightarrow E$ be two genus 2 covers. Then an \emph{isomorphism} between $f_1$ and $f_2$ is an $S$-isomorphism $\varphi : C_1 \longrightarrow C_2$ such that $f_1 = f_2 \circ \varphi$.
\end{definition}

The following lemma shows (in particular) that given two isomorphic genus 2 covers of the same elliptic curve, one of the covers is normalized if and only if the other is.

\begin{lemma}
\label{normalized}
Let $E_1/S, E_2/S$ be an elliptic curves, let $C_1/S, C_2/S$ be genus 2 curves. Let $f : C_2 \longrightarrow E_2$ be a normalized cover, let $\varphi : C_1 \longrightarrow C_2$ be an $S$-isomorphism and $\alpha : E_2 \longrightarrow E_1$ an isomorphism of elliptic curves. Then $\alpha \circ f \circ \varphi: C_1 \longrightarrow E_1$ is normalized.
\end{lemma}
\emph{Proof.}
We can assume that $S$ is connected. Obviously, $\alpha \circ f \circ \varphi$ is minimal. By Corollary \ref{hyerelliptic-involution-commutativity} and (\ref{pseudo-normalized}), we have $\alpha \circ f \circ \varphi \circ \sigma_{C_1/S} = \alpha \circ f \circ \sigma_{C_2/S} \circ \varphi = \alpha \circ [-1]_{E_2/S} \circ f \circ \varphi = [-1]_{E_1/S} \circ \alpha \circ f \circ \varphi : C_1 \longrightarrow E_1$. By \cite[Theorem 3.2 (c)]{KaHur} we have to show that for some geometric point $s \in S$, $(\alpha \circ f \circ \varphi)_s : (C_1)_s \longrightarrow (E_1)_s$ is normalized.\footnote{In \cite[Theorem 3.2 (c)]{KaHur}, the condition that $S$ be connected should be inserted.}

Let $s \in S$. It is well-known that $\varphi_s^{-1}(W_{(C_2)_s}) = W_{(C_1)_s}$. We have $\#(f^{-1}([0_{(E_2)_s}]) \, \cap W_{(C_2)_s}) = \#(\varphi_s^{-1} \left( f^{-1} (\alpha^{-1} ([0_{(E_1)_s}]))  \cap W_{(C_2)_s} \right) ) = \linebreak \# (\varphi_s^{-1} (f^{-1} (\alpha^{-1} ([0_{(E_1)_s}]))) \, \cap \varphi_s^{-1}(W_{(C_2)_s})) = \# ((\alpha \circ f \circ \varphi_s)^{-1} ([0_{(E_1)_s}]) \, \cap W_{(C_1)_s})$. Now with \cite[Corollary 2.3]{KaHur}, the result follows.
\qed\\

The following proposition can be viewed as a Torelli theorem for normalized genus 2 covers of (relative) elliptic curves.
\begin{proposition}
\label{genus-2-cover-Torelli}
Let $E/S$ be an elliptic curve, and let $f_1 : C_1 \longrightarrow E, f_2 : C_2 \longrightarrow E$ be two normalized genus 2 covers. Then the bijection \linebreak $\Iso_S(C_1,C_2) \longrightarrow \Iso_S((J_{C_1},\lambda_{C_1}), (J_{C_2},\lambda_{C_2})), \varphi \mapsto \varphi_*$ of Theorem \ref{genus-2-torelli} induces a bijection between
\begin{itemize}
\item
the set of isomorphisms between the normalized genus 2 covers $f_1$ and \nolinebreak $f_2$

and 
\item
the set of isomorphisms $\alpha$ between the principally polarized abelian varieties $(J_{C_1},\lambda_{C_1})$ and $(J_{C_2},\lambda_{C_2})$ satisfying $(f_1)_* = (f_2)_* \circ \alpha$.
\end{itemize}
\end{proposition}
\emph{Proof.}
We only have to show the surjectivity.

Let $\alpha$ be an isomorphism between $(J_{C_1},\lambda_{C_1})$ and $(J_{C_2},\lambda_{C_2})$ satisfying $(f_1)_* = (f_2)_* \circ \alpha : J_{C_1} \longrightarrow E$. Let $\varphi$ be the unique $S$-isomorphism $C_1 \longrightarrow C_2$ with $\varphi_* = \alpha$. We thus have $(f_1)_* = (f_2 \circ \varphi)_*$. By \cite[Lemma 7.2]{KaHur}, there exists a unique $x \in E(S)$ such that $T_x \circ f_1= f_2 \circ \varphi$. As by Lemma \ref{normalized} both $f_1$ and $f_2 \circ \varphi$ are normalized, we have in fact $f_1 = f_2 \circ \varphi$.
\qed

\begin{remark}
\label{lower-star-upper-star}
The equality $(f_1)_* = (f_2)_* \circ \alpha$ in the above proposition can be restated as $\alpha \circ f_1^* = f_2^*$; cf.\ the calculation in the proof of \cite[Theorem 2.6]{KaHur}.
\end{remark}

\begin{remark}
If $\deg(f_1) \geq 3$ (or $\deg(f_2) \geq 3$), there is in fact at most one isomorphism between $f_1$ and $f_2$; cf.\ \cite[Proposition 3.3]{KaHur}.
\end{remark}

\subsubsection*{Application to the study of the Hurwitz functor}
As in \cite{KaHur}, let $E/K$ be an elliptic curve over a field of characteristic $\neq 2$ (or more generally over a ring in which 2 is invertible or even a scheme over $\mathbb{Z}[1/2]$). As always, we use the notation of \cite{KaHur}.

Proposition \ref{genus-2-cover-Torelli} and Remark \ref{lower-star-upper-star} immediately imply that the ``Torelli map'' $\tau: \mathcal{H}_{E/K,N} \longrightarrow \mathcal{A}_{E/K,N}$ of \cite{KaHur} is a monomorphism for arbitrary $N> 1$; cf.\ \cite[Proposition 5.12]{KaHur}.

The functor $\Psi : \mathcal{H}_{E/K,N} \longrightarrow \mathcal{X}_{E,N,-1}$ of \cite[Corollary 5.13]{KaHur} is thus in fact a monomorphism for arbitrary $N > 1$. Furthermore, the functor $\mathcal{H}_{E/K,N} \longrightarrow \mathcal{J}_{E/K,N}$ of \cite[Proposition 5.17]{KaHur} is an isomorphism for arbitrary $N > 1$, and $\tau : \mathcal{H}_{E/K,N} \longrightarrow \mathcal{A}_{E/K,N}$ is always an open immersion of functors. This of course shortens the proof of Theorem 1.1. at the end of Section 5 in \cite{KaHur}.

It follows that the covers obtained with the ``basic construction'' (\cite[Corollary 5.19]{KaHur}) are always unique up to isomorphism for any $N >1$. For $N \geq 3$, one sees with \cite[Proposition 5.4]{KaHur} that given two covers associated to the same anti-isometry $\psi : E[N] \longrightarrow E'[N]$, there is a \emph{unique} isomorphism between them.

\medskip

Let us again consider genus 2 covers $f: C \longrightarrow E$ of an elliptic curve $E/S$, where $S$ is a scheme over $\mathbb{Z}[1/2]$. The ``basic construction'' now reads as follows (as always, we use the notations of \cite{KaHur}, in particular $E'_f := \ker(f_*)$).

\begin{proposition}[Basic construction]
Let $N> 1$ be a natural number. Let $E/S, E'/S$ be two elliptic curves, and let $\psi : E[N] \longrightarrow E'[N]$ be an anti-isometry which is ``theta-smooth'' (in the sense that the induced principal polarization $\lambda_J$ on $J_\psi := (E \times E')/\Graph(\psi)$ is theta-smooth). Then there is a normalized genus 2 cover $f: C \longrightarrow E$ of degree $N$ such that $(E', \psi)$ is equivalent to $(E'_f,\psi_f)$ (where $\psi_f : E[N] \longrightarrow E'_f[N]$ is the induced anti-isometry). The cover $f$ is unique up to isomorphism (up to unique isomorphism if $N \geq 3$). Moreover, every normalized genus 2 cover of degree $N$ arises in this way.
\end{proposition}

We now give a more symmetric formulation of the ``basic construction''. This ``symmetric basic construction'' has the advantage that it is more rigid than the basic construction for $N=2$.

For this ``symmetric basic construction'', we fix two elliptic curves $E/S$, $E'/S$. 

\begin{definition}
\label{def:symmetric-pair}
A \emph{symmetric pair} (with respect to $E/S$ and $E'/S$) is a triple $(C,f,f')$, where $C/S$ is a genus 2 curve and $f : C \longrightarrow E$, $f' : C \longrightarrow E'$ are minimal covers such that $\ker(f_*) = \Im((f')^*)$ and $\ker(f'_*) = \Im(f^*)$. We say that a symmetric pair is \emph{normalized} if both $f$ and $f'$ are normalized. By an \emph{isomorphism} of two symmetric pairs $(C_1,f_1,f'_1), (C_2, f_2,f'_2)$ we mean an $S$-isomorphism $\varphi : C_1 \longrightarrow C_2$ such that $f_1 = f_2 \circ \varphi$ and $f'_1 = f_2' \circ \varphi$.
\end{definition}

\begin{remark}
It follows from Lemma \ref{normalized} that given two isomorphic symmetric pairs, one of the symmetric pairs is normalized if and only if the other is. 
\end{remark}

\begin{remark}
If $C/S$ is a genus 2 curve and $f : C \longrightarrow E$, $f' : C \longrightarrow E'$ are minimal covers such that $\ker(f_*) = \Im((f')^*)$, then by dualization, one also has $\ker(f'_*) = \Im(f^*)$, i.e.\ $(C,f,f')$ is a symmetric pair.
\end{remark}

\begin{remark}
\label{E'-E'_f-can-iso}
If $(C,f,f')$ is a symmetric pair, then $E'$ (with $(f')^* \circ \lambda_{E'} : E' \longrightarrow J_{C}$) is (canonically isomorphic to) $\ker(f_*) = E'_f$. (If $E/S$ is some elliptic curve, we denote the canonical polarization $E \longrightarrow J_E = \hat{E}$ by $\lambda_E$.)
\end{remark}

\begin{lemma-def}
If $(C,f,f')$ is a symmetric pair, then the degrees of $f$ and $f'$ are equal; this number is called the \emph{degree} of the symmetric pair.
\end{lemma-def}
\emph{Proof.}
Let $N:= \deg(f)$. Then by \cite[Theorem 3.2 (f)]{KaHur}, $f^*$ also has degree $N$. By \cite[Corollary 5.3]{KaHur} and Remark \ref{E'-E'_f-can-iso}, $(f')^* \circ \lambda_{E'} : E' \hookrightarrow J_C$ has also degree $N$, and it follows again with \cite[Theorem 3.2 (f)]{KaHur} that $\deg(f') = \deg((f')^*) = N$.
\qed

\begin{lemma}
\label{c_f}
Let $E/S$ be an elliptic curve, let $C/S$ be a genus 2 curve, and let $f : C \longrightarrow E$ be a minimal cover. Then there exists a unique normalized cover $c_f : C \longrightarrow E'_f$ such that $(c_f)^* \circ \lambda_{E'_f}$ is the canonical immersion $E'_f \hookrightarrow J_C$.\,\footnote{In \cite[Corollary 5.13]{KaHur}, $(c_f)^* \circ \lambda_{E'_f}$ is denoted by $(f')^*$.} In particular, if $f$ is normalized, then $(C,f,c_f)$ is a normalized symmetric pair.
\end{lemma}
\emph{Proof.}
This is a special case of \cite[Theorem 3.2 (f)]{KaHur}.
\qed

\begin{proposition}
\label{symmetric-psi-start}
Let $E/S, E'/S$ be two elliptic curves, let $(C,f,f')$ be a symmetric pair of degree $N$ associated to $E/S$ and $E'/S$. Then there is a unique $\psi : E[N] \tilde{\longrightarrow} E'[N]$ with $(f^{*})_{|E[N]} = (f')^{*} \circ \psi$.\,\footnote{Note that just as in \cite{KaHur} we tacitly identify $E[N]$ with $J_E[N]$.} This $\psi$ is an anti-isometry. Moreover, $\psi$ only depends on the isomorphism class of $(C,f,f')$.
\end{proposition}
\emph{Proof.}
By Remark \ref{E'-E'_f-can-iso}, the existence and uniqueness is \cite[Proposition 5.2]{KaHur}. The fact that $\psi$ only depends on the isomorphism class of $(C,f,f')$ is straightforward.\qed

\begin{proposition}
\label{kernel=Graph(psi)}
With the notation of the previous proposition, let
\[ \pi := f^*  \circ \lambda_{E} \circ \text{\rm pr} + (f')^* \circ \lambda_{E'} \circ \text{\rm pr}' : E \times_S E' \longrightarrow J_C \; ,\]
where $\text{\rm pr} :  E \times_S E' \longrightarrow E$ and $\text{\rm pr} :  E \times_S E' \longrightarrow E'$ are the two projections. Then $\pi$ has kernel $\Graph(-\psi)$. The pull-back to the canonical principal polarization of $J_C$ under $\pi$ is $N$-times the canonical product polarization. In particular, $\psi$ is theta-smooth.
\end{proposition}
\emph{Proof.} This is \cite[Proposition 5.5]{KaHur}.\qed\\

The following ``symmetric basic construction'' can be viewed as a converse to Proposition \ref{symmetric-psi-start}.

\begin{proposition}[Symmetric basic construction]
\label{symmetric-basic-construction}
Let $N>1$ be a \linebreak natural number. Let $E/S, E'/S$ be two elliptic curves, and let $\psi : E[N] \longrightarrow E'[N]$ be an anti-isometry which is theta-smooth. Then there exists a normalized symmetric pair $(C,f,f')$ with respect to $E/S$ and $E'/S$ with $(f^{*})_{|E[N]}$\linebreak$ = (f')^{*} \circ \psi$. The normalized symmetric pair with these properties is essentially unique, i.e.\ it is unique up to unique isomorphism.
\end{proposition}
\emph{Proof.} Let $N$, $E/S,E'/S$ and $\psi : E[N] \longrightarrow E'[N]$ be as in the assertion.

To show the existence, one could use the ``basic construction''. There is however also the following more direct approach:

Consider the abelian variety $J_\psi := (E \times_S E')/\Graph(-\psi)$. By  \cite[Proposition 5.7]{KaHur} there exists a unique principal polarization $\lambda_J$ on $J_\psi$ whose pull-back to $E \times_S E'$ via the projection map is $N$-times the canonical product polarization. By assumption and \cite[Proposition 5.14]{KaHur}, $(J_\psi,\lambda_J)$ is isomorphic to a Jacobian variety of a curve $C/S$. By  \cite[Theorem 3.2 (f)]{KaHur} there exist normalized covers $f : C \longrightarrow E$ and $f' : C \longrightarrow E'$ with $f^*  \circ \lambda_{E} = h_\psi, (f')^* \circ \lambda_{E'} = h_\psi'$, where $h_\psi : E \longrightarrow J_\psi$ and $h_\psi' : E' \longrightarrow J_\psi$ are defined by inclusion into $E \times_S E'$ composed with the projection onto $J_\psi$; cf.\ \cite[Corollary 5.9]{KaHur}. By the exact sequences (28) in \cite[Corollary 5.9]{KaHur}, the conditions $\ker(f_*) = \Im((f')^*)$ and $\ker(f'_*) = \Im(f^*)$ are fulfilled.

We now show the uniqueness. Let $(C_1,f_1,f'_1), (C_2,f_2,f_2')$ be two normalized symmetric pairs associated to $E, E'$ and $\psi$. We claim that there exists a unique isomorphism $\alpha : J_{C_1} \longrightarrow J_{C_2}$ of abelian varieties with $\alpha \circ f_1^* = f_2^*$ and $\alpha \circ (f_1')^* = (f_2')^*$.

Let
\[
\begin{array}{ll}
\pi_1 := f_1^* \circ  \lambda_{E} \circ \text{pr} + (f'_1)^* \circ \lambda_{E'} \circ \text{pr}' : &   E \times_S E' \longrightarrow J_{C_1/S} \; ,\\
\pi_2 := f_2^* \circ \lambda_{E} \circ \text{pr} + (f'_2)^* \circ \lambda_{E'} \circ \text{pr}' : &   E \times_S E' \longrightarrow J_{C_2/S} \; ,
\end{array}
\]
where $\text{pr} : E \times_S E' \longrightarrow E$ and $\text{pr}' : E \times_S E' \longrightarrow E'$ are the two projections.

The two conditions on $\alpha$ are equivalent to $\alpha \circ \pi_1 = \pi_2 : E \times_S E' \longrightarrow J_{C_2/S}$. The assertion follows since by Proposition \ref{kernel=Graph(psi)} $\pi_1 : E \times_S E' \longrightarrow J_{C_1/S} $ and $\pi_2 : E \times_S E' \longrightarrow J_{C_2/S}$ both have kernel $\Graph(-\psi)$.

The fact that $f_1, f_1', f_2$ and $f_2'$ all have degree $N$ implies that the pull-backs of $\lambda_{C_1}$ and $\lambda_{C_2}$ to $E \times_S E'$ via $\pi_1$ and $\pi_2$ respectively are $N$-times the canonical product polarizations. Together with the definition of $\alpha$, this in turn implies that $\hat{\alpha} \circ \lambda_{C_2} \circ \alpha = \lambda_{C_1}$, i.e.\ $\alpha$ preserves the principal polarizations.

Let $\varphi : C_1 \longrightarrow C_2$ be the unique $S$-isomorphism such that $\varphi_* = \alpha$; cf. Theorem \ref{genus-2-torelli}. By Proposition \ref{genus-2-cover-Torelli} and Remark \ref{lower-star-upper-star}, we have $f_1 = f_2 \circ \varphi$ and $f'_1 = f'_2 \circ \varphi$. The uniqueness of $\alpha$ implies that $\varphi : C_1 \longrightarrow C_2$ with these two properties is unique.
\qed

\begin{remark}
\label{largest-theta-smooth}
Let $S$, $E/S$, $E'/S$ and $\psi : E[N] \longrightarrow E'[N]$ be as in the ``symmetric basic construction'' but without the assumption that $\psi$ is theta-smooth. Then by \cite[Corollary 5.16]{KaHur} there exists a uniquely determined largest open subscheme $U$ of $S$ such that $\psi_{|U}$ is theta-smooth. Now $U$ is the largest open subscheme of $S$ over which a symmetric pair with respect to $E_U/U$ and $E'_U/U$ corresponding to $\psi$ exists; this is obvious from Proposition \ref{kernel=Graph(psi)} and the very definition of theta-smoothness.
\end{remark}

\section{Genus 2 covers of degree 2}
\label{section:genus2}
We now concentrate on the case that the covering degree $N$ is $2$. As above, let $S$ be a scheme over $\mathbb{Z}[1/2]$.

In the sequel, by an \emph{isomorphism} $E[2] \longrightarrow E'[2]$, where $E/S$ and $E'/S$ are elliptic curves, we always mean an isomorphism of $S$-group schemes. Note that every such isomorphism is an anti-isogeny. The following proposition is a special case of \cite[Theorem 3]{KaNum}.

\begin{proposition}
\label{theta-smooth-deg-2}
Let $E/S, E'/S$ be two elliptic curves, let $\psi : E[2] \longrightarrow E'[2]$ be an isomorphism. Then $\psi$ is theta-smooth if and only if for no geometric point $s$ of $S$, there exists an isomorphism $\alpha: E_s \longrightarrow E'_s$ such that $\alpha_{|E_s[2]} = \psi_s : E_s[2] \longrightarrow E'_s[2]$.
\end{proposition}

\begin{remark}
Under the conditions of the proposition, let $s$ be a geometric point of $S$. Assume that $E_s$ has $j$-invariant $\neq 0,1728$. Then if $E'_s$ is isomorphic to $E_s$ (i.e.\ if the $j$-invariants of the two curves are equal), there exist exactly two isomorphisms between $E_s$ and $E'_s$. If $\alpha$ is one of these, $-\alpha$ is the other. This means that the isomorphisms between $E_s$ and $E'_s$ induce a \emph{canonical} identification of $E_s[2]$ and $E'_s[2]$. Under the above assumption on the $j$-invariant of $E_s$, the following assertions are thus equivalent.
\begin{itemize}
\item
\emph{There does not exist an isomorphism $\alpha : E_s \longrightarrow E'_s$ such that $\alpha_{|E_s[2]} = \psi_s : E_s[2] \longrightarrow E'_s[2]$.}
\item
\emph{$j(E_s) \neq j(E'_s)$ or $j(E_s) = j(E'_s)$ and, under the canonical identification of $E_s[2]$ and $E'_s[2]$, $\psi_s \neq \id_{E_s[2]}$.}
\end{itemize}
\end{remark}

\begin{proposition}
\label{psi-Weierstrass-deg-2}
Let $E/S, E'/S$ be two elliptic curves with an isomorphism $\psi : E[2] \longrightarrow E'[2]$. Let $C/S$ be a genus 2 curve, and let $(C,f,f')$ be a normalized symmetric pair for $E/S$ and $E'/S$. Then $(f^{*})_{|E[2]} = (f')^{*} \circ \psi$ if and only if $\psi \circ f_{|W_{C/S}} = (f')_{|W_{C/S}}$.
\end{proposition}
\emph{Proof.} Let $E/S, E'/S, \psi, C, f$ and $f'$ be as in the proposition. We only have to show the equivalence after a faithfully flat base change. We can thus assume that $C/S$ has 6 distinct Weierstra{\ss} sections. Now by \cite[Theorem 3.2 (d)]{KaHur}, there exists an embedding $j : C \longrightarrow J_{C}$ which satisfies $j \circ \sigma_{C} = [-1] \circ j$, $[0_{J_{C}}] \cap j(C) = \emptyset$. This implies in particular that $j(W_{C/S}) \subset J_{C}[2]^{\#}$, where $J_C[2]^{\#} := J_C[2] - [0_{J/S}]$.

Assume that $f^*|_{E[2]} = (f')^* \circ \psi$. Then ${f_*}_{|J_{C}[2]} = \lambda_{E}^{-1} \circ (f^*\hat{)} \circ (\lambda_{C})_{|J_{C}[2]} =  \lambda_{E}^{-1} \circ \hat{\psi} \circ ((f')^*\hat{)} \circ (\lambda_{C})_{|J_{C}[2]}$ $= \psi^{-1} \circ {f'_*}_{|J_{C}[2]} : J_C[2] \longrightarrow E[2]$. (We make the usual identification of $E[2]$ with $\hat{E}[2]$ and $J_C[2]$ with $\hat{J}_C[2]$.) Composition with $j_{|W_{C/S}}$ implies $f_{|W_{C/S}} = \psi^{-1} \circ (f')_{|W_{C/S}}$, i.e.\ $\psi \circ f_{|W_{C/S}} = (f')_{|W_{C/S}}$.\\

Let us now assume that $\psi \circ f_{|W_{C/S}} = (f')_{|W_{C/S}}$. We want to show that $\psi \circ {f_*}_{|J_C[2]^{\#}} = {f'_*}_{|J_C[2]^{\#}}$. As $J_C[2] = [0_{J/S}] \stackrel{_\cdot}{\cup} J_C[2]^{\#}$ and clearly $\psi \circ {f_*}_{|[0_{J/S}]} = {f'_*}_{|[0_{J/S}]}$, this implies that $\psi \circ {f_*}_{|J_C[2]} = {f'_*}_{|J_C[2]} : J_C[2] \longrightarrow E[2]$. The equality $(f^*)_{|E[2]} = ((f')^*)_{E[2]} \circ \psi$ then follows by ``dualization'' similarly to above.

By the fact that $(C,f,f')$ is a normalized symmetric pair, we have $\ker(f_*)[2] = \ker(f'_*)[2]$, i.e.\ $\ker({f_*}_{|J_C[2]}) = \ker({f'_*}_{|J_C[2]})$. Let these (equal) kernels be denoted by $K$. Then ${f_*}_{|J_C[2]}$ and ${f'_*}_{|J_C[2]}$ induce homomorphisms $\overline{{f_*}_{|J_C[2]}} : J_C[2]/K \longrightarrow E[2], \overline{{f'_*}_{|J_C[2]}} : J_C[2]/K \longrightarrow E'[2]$. Since these homomorphisms are surjective and $J_C[2]/K$, $E[2]$ and $E'[2]$ are \'etale over $S$ of degree 4, they are in fact isomorphisms. Let 
$p : J_C[2] \longrightarrow J_C[2]/K$ be the canonical projection. Then the equality $\psi \circ {f_*}_{|J[2]} = {f'_*}_{|J[2]}$ implies
\[ \psi \circ \overline{{f_*}_{|J_C[2]}} \circ p \circ j_{|W_{C/S}}= \overline{{f'_*}_{|J_C[2]}} \circ p \circ j_{|W_{C/S}}\; .\]
We claim that $p \circ j_{|W_{C/S}} : W_{C/S} \longrightarrow (J_C[2]/K)^{\#}$ is an \'etale cover.

We have $f_{|W_{C/S}} = \overline{{f_*}_{|J_C[2]}} \circ p \circ j_{|W_{C/S}}$. Since $f_{|W_{C/S}}$ induces an \'etale cover $W_{C/S} \longrightarrow E[2]^{\#}$ of degree 2 and $\overline{{f_*}_{|J_C[2]}}$ is an isomorphism, $p \circ j_{|W_{C/S}} : W_{C/S} \longrightarrow (J_C[2]/K)^{\#}$ is also an \'etale cover of degree 2.

As any surjective \'etale $S$-cover is an epimorphism in the category of $S$-schemes (see \cite[Expos\'e V, Proposition 3.6.]{SGA}), we can thus derive that $\psi \circ \overline{{f_*}_{|J_C[2]}}_{|(J_C[2]/K)^{\#}} = \overline{{f'_*}_{|J_C[2]}}_{|(J_C[2]/K)^{\#}}$, in particular $\psi \circ {f_*}_{|J_C[2]^{\#}} = {f'_*}_{|J_C[2]^{\#}} : J_C[2]^{\#} \longrightarrow E[2]^{\#}$.
\qed\\

With the above two propositions, the ``symmetric basic construction'' can be restated as follows:
\begin{proposition}
\textbf{\emph{(Symmetric basic construction for degree 2 -- \linebreak second form)}}
Let $S$ be a scheme over $\mathbb{Z}[1/2]$. Let $E/S, E'/S$ be two elliptic curves, and let $\psi : E[2] \longrightarrow E'[2]$ be an isomorphism such that for no geometric point $s$ of $S$, there exists an isomorphism $\alpha : E_s \longrightarrow E'_s$ such that $\alpha_{|E_s[2]} = \psi_s$. Then there exists an essentially unique (i.e.\ unique up to unique isomorphism) normalized symmetric pair $(C,f,f')$ with $\psi \circ f_{|W_{C/S}} = (f')_{|W_{C/S}}$.
\end{proposition}

Let $E/S, E'/S$ be elliptic curves, and let $C/S$ be a genus 2 curve. Let $(C,f,f')$ be a normalized symmetric pair with respect to $E/S$ and $E'/S$.

Our goal is now to show that there exists a $\mathbb{P}^1$-bundle $\mathbf{P}$ and covers of degree 2 $E \longrightarrow \mathbf{P}, E' \longrightarrow \mathbf{P}$ such that the induced morphism $C \longrightarrow E \times_{\mathbf{P}} E'$ induces birational morphisms on the fibers over $S$.

Let $\tilde{q} : C \longrightarrow S, q : E \longrightarrow S, q' : E' \longrightarrow S$ be the structure morphisms. Let $\omega_{C/S}:=\tilde{q}_*\Omega_{C/S}$. By Riemann-Roch and ``cohomology and base change'' (\cite[\S 5, Corollary 3]{Mu-AV} and \cite[Theorem 12.11]{Ha}), this is a locally free sheaf of rank 2, and the canonical $S$-morphism $\tilde{\rho} : C \longrightarrow \mathbb{P}(\omega_{C/S})$ is a cover of degree \nolinebreak 2.

By the same general theorems $q_*  \mathcal{L}(2[0_E])$ is a locally free sheaf of rank 2, and the canonical $S$-morphism $\rho : E \longrightarrow \mathbb{P}(q_* \,  \mathcal{L}(2[0_E]) )$ is a cover of degree 2. Analogously, the canonical $S$-morphism $\rho' : E' \longrightarrow \mathbb{P}(q'_* \,  \mathcal{L}(2[0_{E'}]) )$ is a cover of degree 2.

Note that $(C,f,[-1] \circ f'), (C,[-1] \circ f,f')$ and $(C,[-1] \circ f,[-1] \circ f')$ are also normalized symmetric pairs with respect to $E/S$ and $E'/S$ corresponding to $\psi$.

There thus exist unique $S$-automorphisms $\tau, \tau', \tilde{\tau} : C \longrightarrow C$ with 
\[\begin{array}{rclrcl} 
f \circ \tau & = & f         \, , & f' \circ \tau & = & [-1] \circ f' \, ,\\
f \circ \tau' & = & [-1] \circ f \, , & f' \circ \tau' & = & f' \, ,\\
f \circ \tilde{\tau} & = & [-1] \circ f \, , & f' \circ \tilde{\tau} & = & [-1] \circ f' \, .
\end{array}\]
Obviously, $\tau \circ \tau' = \tilde{\tau} = \tau' \circ \tau$ and $\tilde{\tau} = \sigma_{C/S}$.

The automorphisms $\tau$ and $\tau'$ are automorphisms of the covers $f$ and $f'$ respectively, and $\sigma_{C/S}$ is an automorphism of the cover $C \longrightarrow \mathbb{P}(\omega_{C/S})$. We need the following lemma which is a special case of \cite[Lemma 5.6]{LK}.

\begin{lemma}
\label{degree-2-is-quotient}
Let $X$ and $Y$ be connected schemes over $\mathbb{Z}[1/2]$. Let $h : X \longrightarrow Y$ be a finite and flat morphism of degree 2. Then the automorphism group of $h$ is isomorphic to $\mathbb{Z}/2\mathbb{Z}$, and $h$ is a geometric quotient of $X$ under $\Aut(h)$.
\end{lemma}

As a special case of this lemma we obtain: The cover $f : C \longrightarrow E$ is a geometric quotient of $C$ under $\langle \tau \rangle$, and $f' : C \longrightarrow E'$ is a geometric quotient of $C$ under $\langle \tau' \rangle$.

Furthermore, the canonical morphism $\tilde{\rho} : C \longrightarrow \mathbb{P}(\omega_{C/S})$ is a geometric quotient of $C$ under $\langle \sigma_{C/S} \rangle$ (see also \cite[Lemma 3.1]{KaHur} and \cite[Theorem 5.5]{LK}), and the canonical morphisms $\rho : E \longrightarrow \mathbb{P}(q_* \, \mathcal{L}(2[0_E]))$, $\rho' : E' \longrightarrow \mathbb{P}(q'_*\, \mathcal{L}(2[0_{E'}]))$ are geometric quotients of $E$ and $E'$ under $\langle [-1] \rangle$ respectively.

By (\ref{pseudo-normalized}), the automorphism $[-1]$ on $E$ is induced by $\sigma_{C/S}$, and this implies that $\rho \circ f : C \longrightarrow \mathbb{P}(q_* \, \mathcal{L}(2[0_E]))$ is a geometric quotient of $C$ under $\langle \tau, \tau' \rangle = \langle \tau, \sigma_{C/S} \rangle$. Similarly,  $\rho' \circ f' : C \longrightarrow \mathbb{P}(q'_* \, \mathcal{L}(2[0_{E'}]))$ is also a geometric quotient of $C$ under $\langle \tau, \tau' \rangle$. Keeping in mind that a geometric quotient is also a categorial quotient (see \cite[Expos\'e V, Proposition 1.3.]{SGA}), this implies the following theorem.

\begin{theorem}
\label{theorem-P-diagram}
Let $S$ be a scheme over $\mathbb{Z}[1/2]$. Let $C/S$ be a genus 2 curve, $E/S, E'/S$ elliptic curves and $f : C \longrightarrow E, f' : C \longrightarrow E'$ normalized covers of degree 2 with $\ker(f_*) = \Im((f')^*)$, $\ker(f'_*) = \Im(f^*)$. Let $q : C \longrightarrow S, q : E \longrightarrow S, q' : E' \longrightarrow S$ be the structure morphisms, and let $\tilde{\rho} : C \longrightarrow \mathbb{P}(\omega_{C/S})$, $\rho : E \longrightarrow \mathbb{P}(q_* \, \mathcal{L}(2[0_E])), \rho' : E' \longrightarrow \mathbb{P}(q'_* \, \mathcal{L}(2[0_{E'}]))$ be the canonical covers of degree 2.

Then $f$ and $f'$ have unique automorphisms $\tau$ and $\tau'$ respectively which operate non-trivially on all connected components of $C$. These automorphisms have order 2 and satisfy $\tau \circ \tau' = \tau' \circ \tau = \sigma_{C/S}$. The cover $f : C \longrightarrow E$ is a geometric quotient of $C$ under $\langle \tau \rangle$, $f' : C \longrightarrow E'$ is a geometric quotient of $C$ under $\langle \tau' \rangle$, and $\tilde{\rho} : C \longrightarrow  \mathbb{P}(\omega_{C/S})$ is a geometric quotient of $C$ under $\langle \sigma_{C/S} \rangle$.

Now $\rho \circ f: C \longrightarrow {\mathbb{P}(q_* \, \mathcal{L}(2[0_E]))}$ as well as $\rho' \circ f' : C \longrightarrow \mathbb{P}(q'_* \, \mathcal{L}(2[0_{E'}]))$ are geometric quotients of $C$ under $\langle \tau, \tau' \rangle$. We thus have a unique isomorphism $\gamma : {\mathbb{P}(q_* \, \mathcal{L}(2[0_E]))} \longrightarrow {\mathbb{P}(q'_* \, \mathcal{L}(2[0_{E'}]))}$ such that $\gamma \circ \rho \circ f = \rho' \circ f$, and we have unique morphisms $\overline{f} : \mathbb{P}(\omega_{C/S}) \longrightarrow {\mathbb{P}(q_* \, \mathcal{L}(2[0_E]))}$ and $\overline{f'} :  \mathbb{P}(\omega_{C/S}) \longrightarrow {\mathbb{P}(q'_* \, \mathcal{L}(2[0_{E'}]))}$ such that $\rho \circ f = \overline{f} \circ \tilde{\rho}$ and $\rho' \circ f' = \overline{f'} \circ \tilde{\rho}$. All these morphisms are $S$-morphisms, and $\overline{f}, \overline{f'}$ are covers of degree 2.
\end{theorem}
\[\xymatrix{
& C  \ar^{\tilde{\rho}}[d] \ar^{f'}[dr] \ar_f[dl] & \\
E  \ar^{\rho}[d] & {\mathbb{P}(\omega_{C/S})} \ar^{\overline{f'}}[dr] \ar_{\overline{f}}[dl] & {E'} \ar^{\rho'}[d] \\
{\mathbb{P}(q_* \, \mathcal{L}(2[0_E]))} \ar_{\sim}^{\gamma}[rr] & & 
{\mathbb{P}(q'_* \, \mathcal{L}(2[0_{E'}]))}
}\]

\begin{corollary}
\label{corollary-P-diagram}
Let $S$ be a scheme over $\mathbb{Z}[1/2]$, let $C/S$ be a genus 2 curve, let $E/S$ be an elliptic curve, and let $f : C \longrightarrow E$ be a normalized cover of degree 2. Let $\mathbf{P}:=  E/\langle [-1] \rangle = \mathbb{P}(q_* \, \mathcal{L}(2[0_E]))$, let $\rho : E \longrightarrow \mathbf{P}$ be the canonical cover of degree 2, and let $c_f : C \longrightarrow E'_f$ be the normalized cover of degree 2 associated to $f$ by Lemma \ref{c_f}. Then there exists a unique $S$-morphism $\phi' : E'_f \longrightarrow \mathbf{P}$ such that $\rho \circ f = \phi' \circ c_f$. The morphism $\phi'$ is a cover of degree 2.

The induced morphism $C \longrightarrow E \times_{\mathbf{P}} E'_f$ induces birational morphisms on the fibers over $S$.
\end{corollary}

\begin{remark}
Let $S$ be a scheme over $\mathbb{Z}[1/2]$, let $C/S$ be a genus 2 curve, let $E/S$ be an elliptic curve and let $f : C \longrightarrow E$ be a normalized cover of some degree $N$. Let $\tilde{\rho} : C \longrightarrow \mathbb{P}(\omega_{C/S})$, $\rho : E \longrightarrow \mathbb{P}(q_* \, \mathcal{L}(2[0_E]))$ be as above. Then just as in the case that the covering degree is 2, there exists a unique morphism $\overline{f} : \mathbb{P}(\omega_{C/S}) \longrightarrow \mathbb{P}(q_* \, \mathcal{L}(2[0_E]))$ with 
\[\overline{f} \circ \tilde{\rho} = \rho \circ f \; ,\]
and this morphism is a cover of degree $N$.

Indeed, the normalized cover $f$ satisfies $f \circ \sigma_{C/S} = [-1] \circ f$ by (\ref{pseudo-normalized}). This implies that $\rho \circ f \circ \sigma_{C/S} = \rho \circ f$. Note that as above $\tilde{\rho}$ is a geometric quotient of $C$ under $\sigma_{C/S}$. The existence and uniqueness of $\overline{f}$ is now immediate, and it is straightforward to check that $f$ is in fact a cover of degree $N$.
\end{remark}

Let us assume that we are in the situation of the theorem.

The canonical maps $\rho : E \longrightarrow {\mathbb{P}(q_* \, \mathcal{L}(2[0_E]))}$ and $\rho' : E' \longrightarrow {\mathbb{P}(q'_* \, \mathcal{L}(2[0_{E'}]))}$ are ramified at $E[2]$, $E'[2]$ respectively -- these are \'etale covers of $S$ of degree 4 --, and the canonical map $C \longrightarrow {\mathbb{P}(\omega_{C/S})}$ is ramified at $W_{C/S}$ -- this is an \'etale cover of $S$ of degree 6. (We use that $S$ is a scheme over $\mathbb{Z}[1/2]$).

Let $P$ and $P'$ be the relative effective Cartier divisors of $\mathbb{P}(q_* \, \mathcal{L}(2[0_E]))/S$ and $\mathbb{P}(q'_* \, \mathcal{L}(2[0_{E'}]))/S$ associated to the sections $\rho \circ 0_E : S \longrightarrow \mathbb{P}(q_* \, \mathcal{L}(2[0_E]))$ and $\rho' \circ 0_{E'} : S \longrightarrow \mathbb{P}(q'_* \, \mathcal{L}(2[0_{E'}]))$. 

The maps $\rho_{|E[2]^{\#}} : E[2]^{\#} \longrightarrow \mathbb{P}(q_* \, \mathcal{L}(2[0_E]))$ and $(\rho')_{|E'[2]^{\#}} : E'[2]^{\#} \longrightarrow \mathbb{P}(q'_* \, \mathcal{L}(2[0_{E'}]))$ are closed immersions. Let $D$ and $D'$ be the corresponding relative effective Cartier divisors -- they are \'etale covers of degree 3 of $S$. 

Using the theorem, the isomorphism $\psi : E[2] \tilde{\longrightarrow} E'[2]$ corresponding to the isomorphism class of $(C,f,f')$ can be determined in yet another way.

\begin{proposition}
\label{new-psi-description}
Let $\psi : E[2] \tilde{\longrightarrow} E'[2]$. Then $\psi \circ f_{|W_{C/S}} = (f')_{|W_{C/S}}$ if and only if  $\rho' \circ \psi_{|E[2]^{\#}}  = \gamma \circ \rho_{|E[2]^{\#}}$.
\end{proposition}
\emph{Proof.}
The equality $\psi \circ f_{|W_{C/S}} = (f')_{|W_{C/S}}$ implies $\rho' \circ \psi \circ f_{|W_{C/S}} = \rho' \circ (f')_{|W_{C/S}}$, and this implies $\rho' \circ \psi \circ f_{|W_{C/S}} = \gamma \circ \rho \circ f_{|W_{C/S}}$. As $f_{|W_{C/S}} : W_{C/S} \longrightarrow E[2]^{\#}$ is an \'etale cover of degree 2 (thus in particular an epimorphism in the category of \'etale $S$-covers) and $\rho' \circ \psi_{|E[2]^{\#}} : E[2]^{\#} \longrightarrow D'$ as well as $\gamma \circ \rho_{|E[2]^{\#}} : E[2]^{\#} \longrightarrow D'$ are isomorphisms, we can conclude that $\rho' \circ \psi_{|E[2]} = \gamma \circ \rho_{|E[2]^{\#}}$.

Now let $\psi : E[2] \longrightarrow E'[2]$ satisfy $\rho' \circ \psi_{|E[2]^{\#}} = \gamma \circ \rho_{|E[2]^{\#}}$. We have $\rho' \circ \psi \circ f_{|W_{C/S}}= \gamma \circ \rho \circ f_{|W_{C/S}} = \rho' \circ (f')_{|W_{C/S}}$. As $(\rho')_{|E[2]^{\#}} : E'[2]^{\#} \longrightarrow D'$ is an isomorphism, this implies that $\psi \circ f_{|W_{C/S}} = (f')_{|W_{C/S}}$.\qed\\

Let $V$ be the \emph{K\"ahler different divisor} of $f$. By definition, this is the closed subscheme of $C$ which is defined by the zero'th Fitting ideal $F^0(\Omega_{C/E})$ of $\Omega_{C/E} = \Omega_{f}$. (For further information on K\"ahler different divisors see \cite{Kunz}, \cite{LK} or the appendix of \cite{KaGenus2Num}.)

In Section 6 of \cite{LK}, the Weierstra{\ss} divisor of a relative hyperelliptic curve $H/S$ has been defined as the K\"ahler different divisor of the canonical map $H \longrightarrow \mathbb{P}(\omega_{H/S})$. Now the discussion starting at the exact sequence (6.2) until the end of section 6 in \cite{LK} carries over to our case (the only difference being that $V$ has degree 2 and not $2g+2$ over $S$). We thus have:
\begin{lemma} \text{ }
\begin{itemize}
\item
$F^0(\Omega_{C/E}) = \Ann(\Omega_{C/E})$.
\item
$V$ is a relative effective Cartier divisor of degree 2 over $S$.
\item
$V$ is the fixed point subscheme of $C$ under the action of $\tau$, i.e.\ $V$ is the largest subscheme of $C$ with the property that $\tau$ restricts to $V$ and $\tau_{|V} = \id_{V}$.
\item
$V$ is \'etale over $S$.
\end{itemize}
\end{lemma}
\emph{Proof.}
The first assertion, which is written in \cite[Remark 6.4]{LK}, follows from the exact sequence (6.2) in \cite{LK} and the definition of the K\"ahler different divisor. The second, third and forth assertion can be adopted from the text below (6.2) in \cite{LK}, \cite[Proposition 6.5]{LK} and \cite[Proposition 6.8]{LK} respectively.\qed

\begin{lemma}
\label{ram-locus-normal}
If $S$ is reduced, then $V$ is equal to the ramification locus of $f$ endowed with the reduced induced scheme structure.
\end{lemma}
\emph{Proof.}
By the first assertion the previous lemma, the support of $V$ is equal to the set of points where $f$ is ramified, i.e.\ to the ramification locus of $f$. Now since $S$ is reduced and by the previous lemma $V$ is \'etale over $S$, $V$ is reduced (see \cite[Expos\'e I, Proposition 9.2.]{SGA}), and so the assertion follows.
\qed

\begin{proposition}
\label{E-ram-zero}
Under the conditions of Theorem \ref{theorem-P-diagram}, let $\iota : V \hookrightarrow C$ be the canonical closed immersion. Then $(f')_{|V} = f' \circ \iota : V \longrightarrow E'$ is the zero-element in the abelian group $E'(V)$.
\end{proposition}
\emph{Proof.} 
Let $p : V \longrightarrow S$ be the canonical morphism. We have to show that $f' \circ \iota = 0_{E'} \circ p$.

The fact that $\tau_{|V} = \id_{V}$ implies that $[-1] \circ f' \circ \iota = f' \circ \tau \circ \iota = f' \circ \iota$. As $E'[2]$ is the largest closed subscheme $X$ of $E'$ with $[-1]_{|X} = \id_{X}$, this implies that $f' \circ \iota$ factors through $E'[2]$.

Let us now assume that $S$ is connected and let $s$ be some geometric point of $S$. As $E'[2]$ and $V$ are \'etale over $S$, the map $E'[2](V) \longrightarrow E_s'[2](V_s)$ is injective. We thus only have to check that $(f' \circ \iota)_s = 0_{E'_s} \circ p_s : V_s \longrightarrow E'_s$, i.e.\ $f'_s(V_s) = [0_{E'_s}]$. This is equation (\ref{V->0}) in Appendix \ref{over-fields:section}.
\qed

\begin{remark}
Essentially the same statement as in the above proposition holds if $V$ is replaced by the ramification locus endowed with the reduced induced scheme structure (independently of $S$ being reduced). This follows immediately from the proposition because by definition the canonical immersion of this scheme into $C$ factors through $V$.
\end{remark}

\begin{remark}
Let $\Delta := f_*(V)$ be the discriminant divisor of $f$. Then $\Delta$ is a relative effective Cartier divisor of $E/S$ of degree 2. As the geometric fibers over $S$ consist of exactly 2 topological points, it is also \'etale of degree 2 over $S$. In particular, the map $f_{|V} : V \longrightarrow \Delta$ is an isomorphism. Furthermore, if $S$ is reduced, $\Delta$ is equal to the branch locus of $f$ endowed with the reduced induced scheme structure. This can be proved analogously to Lemma \ref{ram-locus-normal}.
\end{remark}

\section{A reformulation of Theorem \ref{theorem-P-diagram}}

Together with the ``symmetric basic construction'' (Proposition \ref{symmetric-basic-construction}) and Proposition \ref{new-psi-description}, a consequence of Theorem \ref{theorem-P-diagram} is:

\emph{
Let $S$ be a scheme over $\mathbb{Z}[1/2]$, and let $E/S, E'/S$ be two elliptic curves and $\psi : E[2] \longrightarrow E'[2]$ a theta-smooth isomorphism. Then with the notations of the previous sections, there is an $S$-isomorphism $\gamma : \mathbb{P}(q_* \, \mathcal{L}(2[0_E])) \tilde{\longrightarrow } $\linebreak${\mathbb{P}(q'_* \, \mathcal{L}(2[0_{E'}]))}$ such that $\rho' \circ \psi_{|E[2]^{\#}} = \gamma \circ \rho_{|E[2]^{\#}}$ holds.
}

The existence of this isomorphism, which is canonically attached to $(E, E', \psi)$ maybe at first sight seems a little bit a mystery. In fact, it can easily be derived from a general statement on $\mathbb{P}^1$-bundles:

Let $E/S, E'/S$ be two elliptic curves with an isomorphism $\psi : E[2] {\longrightarrow} E'[2]$ (not necessarily theta-smooth). Let $\rho : E \longrightarrow \mathbb{P}(q_* \, \mathcal{L}(2[0_E])), \rho' : E' \longrightarrow {\mathbb{P}(q'_* \, \mathcal{L}(2[0_{E'}]))}$ be the corresponding canonical projections. The maps $\rho$ and $\rho'$ are ramified at $E[2]$ and $E'[2]$ respectively. In particular, $\rho_{|E[2]^{\#}} : E[2]^{\#} \hookrightarrow \mathbb{P}(q_* \, \mathcal{L}(2[0_E]))$ and $(\rho')_{|E'[2]^{\#}} : E'[2]^{\#} \hookrightarrow {\mathbb{P}(q'_* \, \mathcal{L}(2[0_{E'}]))}$ are closed immersions. Let $D$ and $D'$ be the corresponding closed subschemes -- these are \'etale covers of $S$ of degree 3. (We use that $S$ is a scheme over $\mathbb{Z}[1/2]$.) Now $\psi_{|E[2]^{\#}} : E[2]^{\#} \tilde{\longrightarrow} E'[2]^{\#}$ induces a canonical isomorphism between $D$ and $D'$. With Proposition \ref{vb-iso}, we conclude:

\begin{proposition}
\label{iso-extension}
There is a unique $S$-isomorphism $\gamma : \mathbb{P}(q_* \, \mathcal{L}(2[0_E])) \tilde{\longrightarrow} \linebreak {\mathbb{P}(q'_* \,  \mathcal{L}(2[0_{E'}]))}$ such that the equality $\rho' \circ \psi_{|E[2]^{\#}} = \gamma \circ \rho_{|E[2]^{\#}}$ holds.
\end{proposition}

Let us again assume that $\psi : E[2] \longrightarrow E'[2]$ is theta-smooth, and let $\gamma$ be as in the proposition. Then we have the following alternative criterion for a triple $(C,f,f')$ to be a normalized symmetric pair.

\begin{proposition}
Let $C/S$ be a genus 2 curve, let $f : C \longrightarrow E, f' : C \longrightarrow E'$ be covers of degree 2. Then $(C,f,f')$ is a normalized symmetric pair corresponding to $\psi$ if and only if $\gamma \circ \rho \circ f = \rho' \circ f'$.
\end{proposition}
\emph{Proof.}
By Theorem \ref{theorem-P-diagram}, Proposition \ref{new-psi-description} and the uniqueness of $\gamma$, it is immediate that a normalized symmetric pair $(C,f,f')$ corresponding to $\gamma$ satisfies $\gamma \circ \rho \circ f = \rho' \circ f' : C \longrightarrow {\mathbb{P}(q'_* \, \mathcal{L}(2[0_{E'}]))}$.

Let this equality be satisfied. If $S$ is the spectrum of an algebraically closed field, the statement is proved in Lemma \ref{over-fields:symmetric-pair}.

In the general case, we can assume that $S$ is connected. As a morphism between (relative) elliptic curves over a connected base is either an isogeny or zero and we already know that $f_* \circ (f')^*$ is zero fiberwise, $f_* \circ (f')^*$ is zero. As $f'$ is obviously minimal, this implies that $\ker(f_*) = \Im((f')^*)$. Similarly, we have $\ker(f'_*) = \Im(f^*)$.

We now want to show that $f$ is normalized. Let $\tau$ be the unique non-trivial automorphism of $f$ which exists by Lemma \ref{degree-2-is-quotient}, similarly let $\tau'$ be the unique non-trivial automorphism of $f'$. Then $\tau \circ \tau' = \sigma_{C/S}, \tau' \circ \tau = \sigma_{C/S}$. (It is not difficult to check these equalities fiberwise, and this suffices by \cite[Lemma 3.1]{KaHur}.)

We claim that $[-1] \circ f = f \circ \sigma_{C/S}$. Indeed, as $\tau \circ \tau' = \tau' \circ \tau$, $\tau'$ induces an automorphism on $E$ over $\mathbb{P}(q_* \, \mathcal{L}(2[0_{E}]))$. By looking at the fibers, one sees that this is not the trivial automorphism. It follows that the induced automorphism is $[-1]$. We thus have $[-1] \circ f = f \circ \sigma_{C/S}$.

By \cite[Theorem 3.2]{KaHur} to show that $f$ is normalized it now suffices to check that for some $s \in S$, $f_s : C_s \longrightarrow E_s$ is normalized. For this statement, we again refer to Lemma \ref{over-fields:symmetric-pair}.

The proof that $f'$ is normalized is analogous.

We have $\rho' \circ \psi \circ f_{|W_{C/S}} = \gamma \circ \rho \circ f_{|W_{C/S}} = \rho' \circ (f')_{|W_{C/S}}$. As $(\rho')_{|E'[2]^{\#}} : E'[2]^{\#} \longrightarrow D$ is an isomorphism, it follows that that $\psi \circ f_{|W_{C/S}} = (f')_{|W_{C/S}}$.

By Proposition \ref{psi-Weierstrass-deg-2}, $(C,f,f')$ is a normalized symmetric pair corresponding to $\psi$.
\qed\\

With the help of Lemma \ref{over-fields:genus-2-deg-2}, we can give a third form of the ``symmetric basic construction'' for $N=2$.

\begin{proposition} 
\textbf{\emph{(Symmetric basic construction for degree 2 -- \linebreak third form)}} 
Let $S$ be a scheme over $\mathbb{Z}[1/2]$, let $E/S$, $E'/S$ be two elliptic curves, and let $\psi : E[2] \longrightarrow E'[2]$ be an isomorphism. Let $\rho : E\longrightarrow \mathbb{P}(q_* \, \mathcal{L}(2[0_E]))$, $\rho' : E' \longrightarrow \mathbb{P}(q'_* \, \mathcal{L}(2[0_{E'}]))$ be the canonical covers of degree 2. Let $\gamma : \mathbb{P}(q_* \, \mathcal{L}(2[0_E])) \longrightarrow \mathbb{P}(q'_* \, \mathcal{L}(2[0_{E'}]))$ be the unique $S$-isomorphism which satisfies $\rho' \circ \psi_{|E[2]^{\#}} = \gamma \circ \rho_{|E[2]^{\#}}$. Assume the following two equivalent conditions are satisfied:
\begin{itemize}
\item
For no geometric point $s$ of $S$, there exists an isomorphism $\alpha : E_s \longrightarrow E'_s$ with $\alpha_{|E_s[2]} = \psi_s$.
\item
The images of the sections $\rho' \circ 0_{E'}$ and $\gamma \circ \rho \circ 0_E$ of $\mathbb{P}(q'_* \, \mathcal{L}(2 [0_{E'}] )) \longrightarrow S$ are disjoint.
\end{itemize}
Then there exists a curve $C/S$ and covers $f : C \longrightarrow E$, $f' : C \longrightarrow E'$ of degree 2 such that $\gamma \circ \rho \circ f = f' \circ \rho'$. Any such triple $(C,f,f')$ is a normalized symmetric pair corresponding to $\psi$, and it is unique up to unique isomorphism.
\end{proposition}

If one assumes that the base-scheme is regular, one can give a more concrete description of the curve $C$ and the covers $f,f'$ (as well as to prove its existence in an alternative way).

\begin{proposition}
Under the conditions of the above proposition, let $S$ be regular. Then $E \times_{\mathbb{P}(q'_* \, \mathcal{L}(2[0_{E'}]))} E'$ (where the product is with respect to $\gamma \circ \rho$ and $\rho'$) is reduced with total quotient ring $\kappa(E) \times_{\kappa(\mathbb{P}(q'_* \, \mathcal{L}(2[0_{E'}])))} \kappa(E')$. The normalization $C$ of $E \times_{\mathbb{P}(q'_* \, \mathcal{L}(2[0_{E'}]))} E'$ is a genus 2 curve, and the induced maps $f: C \longrightarrow E, f' : C \longrightarrow E'$ are degree 2 covers which satisfy  $\gamma \circ \rho \circ f = f' \circ \rho'$.
\end{proposition}
\emph{Proof.}
As $S$ is regular, it is also locally integral, in particular, its connected components are integral; see \cite[Theorem 14.3]{Ma},  \cite[I (4.5.6)]{EGA}. We can thus assume that $S$ is integral.

Let $\mathcal{F} := \rho'_* \, \mathcal{L}(2[0_{E'}])$. We first show that $E \times_{\mathbb{P}(\mathcal{F})} E'$ is integral and that its function field is $\kappa(E) \otimes_{\kappa(\mathbb{P}(\mathcal{F}))} \kappa(E')$.

The ring $\kappa(E) \otimes_{\kappa(\mathbb{P}(\mathcal{F}))} \kappa(E')$ is a field because by assumption, the generic points of $\rho'([0_{E'}])$ and $\gamma(\rho([0_E]))$ are distinct.

Let $A$ be the coordinate ring of an affine open part $U$ of $\mathbb{P}(\mathcal{F})$, let $B$ and $C$ the corresponding rings of the preimages of $U$ in $E$ and $E'$. We claim that the canonical map $B \otimes_A C \longrightarrow \kappa(B) \otimes_{\kappa(A)} \kappa(C) \simeq \kappa(E) \times_{\kappa(\mathbb{P}(\mathcal{F}))} \kappa(E')$ is injective.

We have $\kappa(B) \otimes_{\kappa(A)} \kappa(C) \simeq (B \otimes_A C) \otimes_A \kappa(A)$ as $B$ and $C$ are finite over $A$. We thus have to show that the map $A \otimes_B C \longrightarrow (B \otimes_A C) \otimes_A \kappa(A)$ is injective. Now, $A \longrightarrow \kappa(A)$ is injective and $B \otimes_A C$ is flat over $A$ ($C$ is flat over $A$, thus $C \otimes_A B$ is flat over $B$, and as $B$ is flat over $A$, $B \otimes_A C$ is flat over $A$). This implies that $B \otimes_A C \longrightarrow (B \otimes_A C) \otimes_A \kappa(A)$ is injective. It follows that $B \otimes_A C$ is reduced.

We have seen that $B \otimes_A C$ is contained in the field $(B \otimes_A C) \otimes_A \kappa(A)$, and obviously $(B \otimes_A C) \otimes_A \kappa(A)$ is contained in the function field of $B \otimes_A C$. This implies that $(B \otimes_A C) \otimes_A \kappa(A) \simeq \kappa(E) \otimes_{\kappa(\mathbb{P}(\mathcal{F}))} \kappa(E')$ is the function field of $B \otimes_A C$. 

We have seen that $E \times_{\mathbb{P}(\mathcal{F})} E'$ is integral (in particular reduced) and its function field is indeed $\kappa(E) \otimes_{\kappa(\mathbb{P}(\mathcal{F}))} \kappa(E')$.

\smallskip

We now show the statements on $C$.

The field $\kappa(S)$ is algebraically closed in $\kappa(E) \times_{\mathbb{P}(\kappa(\mathcal{F}))} \kappa(E')$, and as $S$ is regular, $S$ is normal; see \cite[Theorem 19.4]{Ma}. This implies with \cite[III (4.3.12)]{EGA} that the geometric fibers of $C$ over $S$ are connected.

Let $W$ be the different divisor of $E \times_{\mathbb{P}(\mathcal{F})} E' \longrightarrow \mathbb{P}(\mathcal{F})$. Then $(E \times_{\mathbb{P}^1} E') - W$ is normal, because the domain of an \'etale morphism mapping to a normal scheme is normal; see \cite[Expos\'e I, Corollaire 9.11.]{SGA}. It follows that $C \longrightarrow E \times_{\mathbb{P}(\mathcal{F})} E'$ induces an isomorphism between the complement of the preimage of $W$ in $C$ and $(E \times_{\mathbb{P}(\mathcal{F})} E') - W$. Since the restriction of $W$ to the fibers over $S$ is zero-dimensional, it follows that $C \longrightarrow E \times_{\mathbb{P}(\mathcal{F})} E'$ induces birational morphisms on the fibers over $S$.

By Abhyankar's Lemma (\cite[Expos\'{e} X, Lemme 3.6]{SGA}) and ``purity of the branch locus'' (\cite[Expos\'{e} X, Th\'eor\`eme 3.1.]{SGA}), $f$ is \'etale outside $(f')^{-1}([0_{E'}]))$ and $f'$ is \'etale outside $f^{-1}([0_{E}])$. Let $x$ be a topological point of $C$. As by assumption $\gamma(\rho([0_E]))$ and $\rho'([0_{E'}])$ are disjoint, $x \notin (f')^{-1}([0_{E'}])$ or $x \notin f^{-1}([0_E])$. In the first case, the morphism $f$ is \'etale at $x$, and since $E$ is smooth over $S$, $C$ over $S$ is smooth at $x$. In the second case, the argument is analogous and the conclusion is the same. It follows that $C$ is smooth over $S$.

Let $s$ be a geometric point of $S$. We have already shown that $C_s$ is connected, and by what we have just seen, $C_s$ is non-singular. We have to show that the genus of this curve is 2. We already know that $C_s \longrightarrow E_s \times_{\mathbb{P}^1_{\kappa(s)}} E'_{s}$ is birational. It follows that $C_s \longrightarrow E_s$ has degree 2. Since $\gamma(\rho([0_E])) \neq \rho'([0_{E'}])$, the morphism $C_s \longrightarrow E_s$ is ramified exactly at the preimages of $\rho'([0_{E'}])$ in $E_s$ (here we use again Abhyankar's Lemma). This preimage consists of exactly two closed points. It follows that the genus of $C_s$ is 2.
\qed

\bigskip

\appendix

\section{Genus 2 covers of degree 2 over fields}
\label{over-fields:section}
In this part of the appendix, we provide some results on genus 2 covers of elliptic curves of degree 2 over algebraically closed fields of characteristic $\neq 2$.

In the following, let $\overline{\kappa}$ be an algebraically closed field of characteristic $\neq 2$. Let $E/\overline{\kappa}, E'/\overline{\kappa}$ be two elliptic curves, $\psi : E[2] \tilde{\longrightarrow} E'[2]$. Let $\phi : E \longrightarrow \mathbb{P}^1_{\overline{\kappa}}, \phi' : E' \longrightarrow \mathbb{P}^1_{\overline{\kappa}}$ be two covers of degree 2 which are ramified at $E[2]$ and $E'[2]$ respectively such that $\phi' \circ \psi_{|E[2]^{\#}} = \phi_{|E[2]^{\#}}$. Let $C$ be the normalization of $E \times_{\mathbb{P}^1_{\overline{\kappa}}} E'$.

Let $P:= \phi([0_E])$, $P' := \phi'([0_{E'}])$. By assumption, $\rho(E[2]^{\#}) = \rho'(E'[2]^{\#})$; let this divisor be denoted by $D$.

\begin{lemma}
\label{over-fields:genus-2-deg-2}
The following assertions are equivalent.
\begin{enumerate}[a)]
\item
The points $P$ and $P'$ are distinct.
\item
$E \times_{\mathbb{P}^1_{\overline{\kappa}}} E'$ is irreducible.
\item
$C/\overline{\kappa}$ is a genus 2 curve.
\item
The two covers $\phi: E \longrightarrow \mathbb{P}^1_{\overline{\kappa}}$ and $\phi' : E' \longrightarrow \mathbb{P}^1_{\overline{\kappa}}$ are not isomorphic (i.e.\ there does not exist a $\overline{\kappa}$-isomorphism $\alpha : E \longrightarrow E'$ with $\phi = \phi' \circ \alpha$).
\item
There does not exist an isomorphism of elliptic curves $\alpha: E \longrightarrow E'$ with $\alpha_{|E[2]} = \psi$
.
\end{enumerate}
\end{lemma}
\emph{Proof.} Keeping in mind that $C$ is regular, i.e.\ smooth over $\Spec(\overline{\kappa})$, the equivalence of the first four assertions is not difficult to show.

Assume that the covers are isomorphic via $\alpha : E \longrightarrow E'$. Then in particular $P=P'$.  We have the isomorphisms $\phi_{|E[2]} : E[2] \longrightarrow D \cup P$, $(\phi')_{|E[2]} : E[2] \longrightarrow D \cup P$. It follows that $\alpha_{|E[2]} = ({\phi'_{|D \cup P}})^{-1} \circ \phi_{E[2]} = \psi$. In particular, $\alpha$ is an isomorphism of elliptic curves.

On the other hand, assume that there exists an isomorphism of elliptic curves $\alpha : E \longrightarrow E'$ with $\alpha_{|E[2]} = \psi$. Then $\phi_{|E[2]} = \phi' \circ \alpha_{|E[2]}$. It is well-known that this implies that $\phi = \phi' \circ \alpha$.
\qed\\

Let us assume that the equivalent conditions of the lemma are satisfied. Then we have a commutative diagram
\begin{equation}
\label{over-fields:diagram}
\xymatrix{
& \ar_f[dl] C \ar_{\tilde{\phi}}[d] \ar^{f'}[dr] & \\
E \ar_{\phi}[dr] & \ar_{\overline{f}}[d] {\mathbb{P}^1_{\overline{\kappa}}} & {E'\; ,} \ar^{\phi'}[dl] \\
& {\mathbb{P}^1_{\overline{\kappa}}}
}
\end{equation}
where all morphisms are covers of degree 2. We have that
\begin{itemize}
\item
$\overline{f} : \mathbb{P}^1_{\overline{\kappa}} \longrightarrow \mathbb{P}^1_{\overline{\kappa}}$ is branched exactly at the set $P \cup P'$,
\item
$\tilde{\phi} : C \longrightarrow \mathbb{P}^1_{\overline{\kappa}}$ is branched exactly at the set $\overline{f}^{-1}(D)$,
\item
$f : C \longrightarrow E$ is branched exactly at the set $\phi^{-1}(P')$,
\item
$f' : C \longrightarrow E'$ is branched exactly at the set $(\phi')^{-1}(P')$.
\end{itemize}
These statements can for example easily be proved with Abhyankar's Lemma.

Let $V \subset C$ be the ramification locus of $f$. Then $(\phi \circ f)(V) = P'$, i.e.\ $(\phi' \circ f')(V) = P'$, and this implies
\begin{equation}
\label{V->0}
f'(V) = [0_{E'}] \; .
\end{equation}

\begin{lemma}
\label{over-fields:symmetric-pair}
$(C,f,f')$ is a normalized symmetric pair with respect to $E$ and $E'$ corresponding to $\psi$.
\end{lemma}
\emph{Proof.}
It is not difficult to show that we have a commutative diagram
\[ \xymatrix{ 
& \ar_{f_*}[dl] {J_{C_s}} & \\
{J_{E}} && {J_{E'}} \ar_{(f')^*}[ul] \ar^{(\phi')_*}[dl] \; . \\
& \ar^{\phi^*}[ul] {J_{\mathbb{P}^1_{\kappa}} = 0}}
 \]
This implies that $f_* \circ (f')^*$ is zero. As $f'$ is obviously minimal, this implies that $\ker(f_*) = \Im((f')^*)$. Similarly, we have $\ker(f'_*) = \Im(f^*)$.

By the above statements on the branching of $\overline{f}$ and $\tilde{\phi}$, over each point of $D$, there lie exactly 2 Weierstra{\ss} points. This implies that over each point of $E[2]^{\#}$ there also lie exactly 2 Weierstra{\ss} points. It follows that $f$ is normalized.

The proof that $f'$ is normalized is analogous.

We have $\phi' \circ \psi \circ f_{|W_{C/S}} = \phi \circ f_{|W_{C/S}} = \phi' \circ (f')_{|W_{C/S}}$. As $(\phi')_{|E[2]^{\#}} : E'[2]^{\#} \longrightarrow D$ is an isomorphism, we can conclude that $\psi \circ f_{|W_{C/S}} = (f')_{|W_{C/S}}$.

By Proposition \ref{psi-Weierstrass-deg-2}, it follows that $(C,f,f')$ is a normalized symmetric pair corresponding to $\psi$.
\qed

\begin{remark}
By Proposition \ref{theta-smooth-deg-2}, the last assertion of Lemma \ref{over-fields:genus-2-deg-2} is equivalent to $\psi$ being irreducible (i.e.\ theta-smooth).

Lemmata \ref{over-fields:genus-2-deg-2} and \ref{over-fields:symmetric-pair} can however also be used to prove Proposition \ref{theta-smooth-deg-2} (i.e.\ \cite[Theorem 3]{KaNum} in the special case that the covering degree is 2). By the definition of Theta-smoothness, we can thereby restrict ourselves to 
the case that $S=\overline{k}$.

If $\psi$ satisfies the conditions of Lemma \ref{over-fields:genus-2-deg-2}, then by Lemma \ref{over-fields:symmetric-pair} and Proposition \ref{kernel=Graph(psi)}, $\psi$ is irreducible.

On the other hand, if $\psi$ is irreducible and $(C,f,f')$ is the corresponding symmetric pair, then we have degree 2 covers $\phi : E \longrightarrow \mathbb{P}^1_{\overline{\kappa}}, E' \longrightarrow \mathbb{P}^1_{\overline{\kappa}}$ which ramify at $E[2]$ and $E'[2]$ respectively with $\phi \circ f = \phi' \circ f'$ (for example by Theorem \ref{theorem-P-diagram}). Consequently, the equivalent conditions of Lemma \ref{over-fields:genus-2-deg-2} hold.

Also Remark \ref{largest-theta-smooth} can -- for covering degree 2 -- be derived from Lemma \ref{over-fields:symmetric-pair}: The open subset $U$ of $S$ where $P$ and $P'$ do not meet obviously has the correct properties.
\end{remark}

\section{Some results on projective space bundles}
\label{projective-bundles}
In the following, let $S$ be an arbitrary (not necessarily locally noetherian) scheme. Let $\mathbb{P}^1_S := \text{Proj}(\mathbb{Z}[X_0,X_1]) \times_{\Spec(\mathbb{Z})} S$. Then $\mathcal{O}(1)$ on $\mathbb{P}^1_S$ has two canonical global generators, $X_0$ and $X_1$.
\begin{lemma}
\label{Aut-P^1-3-transitive}
Let $s_1, s_2, s_3, s'_1, s'_2,s'_3 : S \longrightarrow \mathbb{P}^1_S$ be six sections of $\mathbb{P}^1_S \longrightarrow S$ such that the images of $s_1, s_2, s_3$ as well as of $s_1', s_2', s_3'$ are pairwise disjoint. Then there exists a unique $S$-automorphism $\beta$ of $\mathbb{P}^1_S$ with $\beta \circ s_i = s'_i$ for $i=1,2,3$.
\end{lemma}
\emph{Proof.}
By considering an open affine covering, we can restrict ourselves to the case that $S$ is affine. The general case then follows by the uniqueness of $\alpha$.

Each of the $s_i, s'_i$ is given by an invertible sheaf with two global sections which generate it; cf.\ \cite[II, Theorem 7.1.]{Ha}.
Let $U = \Spec(A)$ be an affine open subset such that all these sheaves are trivial. We are going to show the result for $(s_i)_{|U}, (s'_i)_{|U}$ over $U$. Again the result in the lemma then follows by the uniqueness of $\alpha$ on $U$ via the consideration of an open affine covering. Let us denote $(s_i)_{|U}$ by $s_i$, $(s'_i)_{|U}$ by $s'_i$.

If $\beta : \mathbb{P}^1_A \longrightarrow \mathbb{P}^1_A$ is an automorphism, then $\beta^* (\mathcal{O}(1)) \approx \mathcal{O}(1) \otimes p^*(\mathcal{L})$, where $p : \mathbb{P}^1_A \longrightarrow \Spec(A)$ is the structure morphism and $\mathcal{L}$ is an invertible sheaf on $\Spec(A)$; see \cite[0. \S 5 b)]{Mu-GIT}.

Let us assume that $\beta \in \Aut_A(\mathbb{P}^1_A)$ satisfies $\beta \circ s_i = s_i'$ for some $i$, and let $\mathcal{L}$ be as above. Then $\mathcal{L} = (s_i)^* p^*(\mathcal{L}) = (s_i)^* \beta^* (\mathcal{O}(1)) = (s'_i)^*(\mathcal{O}(1)) = \mathcal{O}_{\text{Spec}(A)}$ by the above assumption on $A$.

We can thus restrict ourselves to automorphisms $\beta$ with $\beta^* (\mathcal{O}(1)) \approx \mathcal{O}(1)$. Fixing an isomorphism of $\beta^* (\mathcal{O}(1))$ with $\mathcal{O}(1)$, $\beta^{*}X_0$ and $\beta^{*}X_1$ define two global sections of $\mathcal{O}(1)$. Thus $\beta$ corresponds to two global section of $\mathcal{O}(1)$ which are unique up to multiplication by an element of $A^*$. Such elements can be written as $a X_0 + b X_1, c X_0 + d X_1$ $(a,b,c,d \in A)$ such that the matrix $\left( \begin{array}{cc} a & c \\ b & d \end{array} \right)$ is invertible. The matrix $\left( \begin{array}{cc} a & c \\ b & d \end{array} \right)$ is thereby unique up to multiplication by an element of $A^*$.

By assumption on $U$, any of the sections $s_i, s_i'$ is given by a tuple of two elements of $A$ which generate the unit ideal. Furthermore, each of these tuples is unique up to multiplication by an element of $A^*$. We can thus uniquely represent any of the $s_i, s_i'$ by an element in $A^2/A^*$.

Let $(f,g) \in A^2/A^*$ be such an element corresponding to $s_i$. Then $\beta \circ s_i$ is given by $(fa + gb, fc + gd) \in A^2/A^*$, i.e.\ it is given by the usual application of $\left( \begin{array}{cc} a & c \\ b & d \end{array} \right)$ on $(f,g)$ from the right.

Note that the assumption on the images of the $s_i$ and $s_i'$ is equivalent to the condition that for all $t \in S$, the restrictions of $s_1, s_2, s_3$ to the fiber over $t$ as well as the restrictions of the $s_1',s_2',s_3'$ are distinct. This in turn is equivalent to the condition that for all prime ideals $P$ of $A$, the tuples $(f,g)$ as above stay distinct in $(A/P)^2/(A/P)^*$.

Now the result of this lemma follows from the following lemma which - for convenience - we formulate with the usual left operation. \qed
\\

We introduce the following notation: For $v \in A^2$, we write $\tilde{v}$ for the reduction of $v$ modulo $A^*$.

\begin{lemma}
Let $\left(\begin{array}{c} a_i \\ b_i\end{array}\right), \left(\begin{array}{c} a'_i \\ b'_i \end{array}\right) \in A^2$ for $i=1,2,3$ be given such that for all prime ideals $P$ of $A$, the $\left(\begin{array}{c} a_i \\ b_i \end{array}\right)$ for $i=1,2,3$ as well as the $\left(\begin{array}{c} a'_i \\ b'_i \end{array}\right)$ for $i=1,2,3$ define pairwise distinct elements in $(A/P)^2/(A/P)^*$. Then there exists an invertible matrix $B \in M_{2 \times 2}(A)$, unique up to multiplication by an element of $A^*$, such that $B \widetilde{\left( \begin{array}{c} a_i \\ b_i \end{array} \right)} = \widetilde{\left( \begin{array}{c} a'_i \\ b'_i \end{array} \right)} \in A^2/A^*$.
\end{lemma}
\emph{Proof.}
We show the existence first.

We only have to show the existence for $\widetilde{\left(\begin{array}{c} a_1' \\ b_1' \end{array}\right)} = \widetilde{\left(\begin{array}{c} 1 \\ 0 \end{array}\right)}, \widetilde{\left(\begin{array}{c} a_2' \\ b_2' \end{array}\right)} = \widetilde{\left(\begin{array}{c} 0 \\ 1 \end{array}\right)}, \widetilde{\left(\begin{array}{c} a_3' \\ a_3' \end{array}\right)} = \widetilde{\left(\begin{array}{c} 1 \\ 1 \end{array}\right)}$.

We claim that the matrix $M:= \left( \begin{array}{cc} a_1 & a_2 \\ b_1 & b_2 \end{array} \right) \in M_{2 \times 2}(A)$ is invertible. Let $d$ be the determinant of this matrix. By assumption, for all prime ideals $P$ of $A$, the reduction of $d$ modulo $P$ is non-zero. It follows that $d$ does not lie in any prime ideal, thus it is a unit (as otherwise it would lie in a maximal ideal).

Now $M^{-1}$ maps $\left(\begin{array}{c} a_1 \\ b_1 \end{array} \right)$ to $\left(\begin{array}{c} 1 \\ 0 \end{array}\right)$ and $\left(\begin{array}{c} a_2 \\ b_2 \end{array}\right)$ to $\left(\begin{array}{c} 0 \\ 1 \end{array}\right)$. Let $\left(\begin{array}{c} a \\ b \end{array}\right)$ be the image of $\left(\begin{array}{c} a_3 \\ b_3 \end{array}\right)$. The assumption remains valid for the images of $\left( \begin{array}{c} a_i \\ b_i \end{array} \right)$ under $M^{-1}$, and it says that $a$ and $b$ are not divisible by any prime ideal, i.e.\ they are units. The invertible matrix $M' := \left( \begin{array}{cc} a^{-1} & 0\\ 0 & b^{-1} \end{array} \right)$ fixes $\widetilde{\left(\begin{array}{c} 1 \\ 0 \end{array}\right)}$ and $\widetilde{\left(\begin{array}{c} 0 \\ 1 \end{array}\right)}$ and maps $\left(\begin{array}{c} a \\ b \end{array}\right)$ to $\left(\begin{array}{c} 1 \\ 1 \end{array}\right)$, so $B := M' M^{-1}$ has the desired properties.

Given what we have already shown, for the uniqueness it suffices to remark that only matrixes of the form $aI$ $(a \in A^*)$ fix $\widetilde{\left(\begin{array}{c} 1 \\ 0 \end{array}\right)}, \widetilde{\left(\begin{array}{c} 0 \\ 1 \end{array}\right)}$ and $\widetilde{\left(\begin{array}{c} 1 \\ 1 \end{array}\right)}$.
\qed

\begin{lemma}
Let $D, D'$ be two subschemes of $\mathbb{P}^1_S$ such that $D \longrightarrow S, D' \longrightarrow S$ are \'etale covers of degree 3, let $\eta : D \longrightarrow D'$ be an $S$-isomorphism. Then there exists a unique $S$-automorphism of $\mathbb{P}^1_S$ such that $\alpha_{|D} = \eta$.
\end{lemma}
\emph{Proof.}
As $D \longrightarrow S$ is an \'etale cover, there exists a Galois cover $T \longrightarrow S$ such that $D_T = D \times_S T \simeq T \cup T \cup T$ (isomorphism over $T$); cf.\ \cite[Expos\'e V, 4 g)]{SGA}.

Let $t_1, t_2, t_3 : T \longrightarrow D_T$ be the three immersions. Then for any $\alpha \in \mathbb{P}^1_T$, the condition $\alpha_{|D_T} = \eta_T$ is equivalent to $\alpha \circ t_i = \eta_T \circ t_i$ for $i=1,2,3$.

It follows from Lemma \ref{Aut-P^1-3-transitive} that there exists a unique automorphism $\alpha$ of $\mathbb{P}^1_T$ such that $\alpha_{|D_T} = \eta_T$.

This implies by Galois descent that there exists a unique automorphism $\alpha$ of $\mathbb{P}^1_S$ with $\alpha_{|D} = \eta$.
\qed

\begin{proposition}
\label{vb-iso}
Let $\mathbf{P}, \mathbf{P'}$ be two $\mathbb{P}^1$-bundles over $S$. Let $D$ be a subscheme of $\mathbf{P}$, $D'$ a subscheme of $\mathbf{P'}$ such that $D \longrightarrow S$ and $D' \longrightarrow S$ are \'etale covers of degree 3. Let $\eta : D \longrightarrow D'$ be an $S$-isomorphism. Then there exists a unique $S$-isomorphism $\alpha : \mathbf{P} \longrightarrow \mathbf{P'}$ such that $\alpha_{|D} = \eta$.

In particular, if $\mathbf{P}$ has three sections over $S$ which do not meet, it is $S$-isomorphic to $\mathbb{P}^1_S$.
\end{proposition}
\emph{Proof.}
If $\mathbf{P}$ and $\mathbf{P'}$ are trivial bundles (i.e.\ $S$-isomorphic to $\mathbb{P}^1_S$), the result follows immediately from the previous lemma. The general case follows from the uniqueness of $\alpha$ by a glueing argument.
\qed

\begin{remark}
The subscheme $D$ of $\mathbf{P}$ in the proposition is in fact a relative effective Cartier divisor of $\mathbf{P}$. This follows from \cite[Corollary 3.9]{Mi-JV}.
\end{remark}

\bibliography{basics-genus2-split}

\begin{thebibliography}{10}

\bibitem{BLR}
S.~Bosch, W.~L\"utkebohmert, and M.~Raynaud.
\newblock {\em {N\'eron Models}}.
\newblock Springer-Verlag, 1980.

\bibitem{DM}
P.~Deligne and D.~Mumford.
\newblock The irreducibility of the space of curves of given genus.
\newblock {\em Inst.\ Hautes Etudes Sci.\ Publ.\ Math.}, pages 75--109, 1969.

\bibitem{DF}
C.~Diem and G.~Frey.
\newblock {Non-constant genus 2 curves with pro-Galois covers}.
\newblock To appear.

\bibitem{FKBasicC}
G.~Frey and E.~Kani.
\newblock Curves of genus 2 covering elliptic curves and an arithmetical
  application.
\newblock In {\em Arithmetic Algebraic Geometry (Texel Conf., 1989)}, volume~89
  of {\em {Prog.\ Math.}}, pages 153--176, Boston, 1991. {Birkh\"auser}.

\bibitem{SGA}
A.~Grothendieck.
\newblock {\em Rev\^etements \'etales et groupe fondamental (SGA I)}, volume
  1960/61 of {\em S\'eminaire de G\'eom\'etrie Alg\'ebrique}.
\newblock Institut des Hautes \'Etudes Scientifiques, Paris.

\bibitem{EGA}
A.~{Grothendieck with J. Dieudonn\'e}.
\newblock {El\'ements de G\'eom\'etrie Al\-g\'e\-brique (I-IV)}.
\newblock {ch.\ I: Springer-Verlag, Berlin, 1971; ch. II-IV: Publ.\ Math.\
  Inst.\ Hautes Etud.\ Sci. 8,11, 17, 20,24, 28, 32, 1961-68}.

\bibitem{Ha}
R.~Hartshorne.
\newblock {\em Algebraic Geometry}.
\newblock Springer-Verlag, New York, 1977.

\bibitem{KaGenus2Num}
E.\ Kani.
\newblock The number of genus 2 covers of an elliptic curve.
\newblock To appear.

\bibitem{KaNum}
E.~Kani.
\newblock The number of curves of genus two with elliptic differentials.
\newblock {\em J.\ reine angew.\ Math.}, 485:93--121, 1997.

\bibitem{KaHur}
E.~Kani.
\newblock Hurwitz spaces of genus 2 covers of an elliptic curve.
\newblock {\em Collect.\ Math.}, 54(1):1--51, 2003.

\bibitem{Kra}
A.~Krazer.
\newblock {\em Lehrbuch der Thetafunktionen}.
\newblock Chelsea Publishing Company, 1970.
\newblock Reprint of the 1903 edition.

\bibitem{Kuhn}
R.~Kuhn.
\newblock {Curves of genus 2 with split Jacobian}.
\newblock {\em Trans.\ Am.\ Math.\ Soc.}, 307:41--49, 1988.

\bibitem{Kunz}
E.\ Kunz.
\newblock {\em {K\"ahler Differentials}}.
\newblock Vieweg, Wiesbaden, 1986.

\bibitem{LK}
K.~L{\o}nsted and S.~Kleiman.
\newblock Basics on families of hyperelliptic curves.
\newblock {\em Compos.\ Math.}, 38:83--111, 1979.

\bibitem{Ma}
H.~Matsumura.
\newblock {\em Commutative ring theory}.
\newblock Cambidge University Press, Cambridge, UK, 1986.

\bibitem{Mi-JV}
J.~Milne.
\newblock Jacobian varieties.
\newblock In G.~Cornell and J.~Silverman, editors, {\em Arithmetic Geometry},
  pages 167--212. Springer-Verlag, 1986.

\bibitem{Mu-GIT}
D.~Mumford.
\newblock {\em {Geometric Invariant Theory}}.
\newblock Springer-Verlag, 1965.

\bibitem{Mu-AV}
D.~Mumford.
\newblock {\em {Abelian Varieties}}.
\newblock Tata Institute for Fundamental Research, 1970.

\end{thebibliography}
\bibliographystyle{plain}

\vspace{4 ex}
\noindent
Universit\"at Leipzig, Fakult\"at f\"ur Mathematik und Informatik, \\
Augustusplatz~10, 04109 Leipzig, Germany. \\
email: diem@math.uni-leipzig.de

\end{document}